\documentclass[twocolumn,10pt]{scrartcl}

\usepackage{amsmath,amstext,amsbsy,amsopn,amscd,amsxtra,upref,amssymb}
\usepackage{graphicx}
\usepackage{dcolumn}
\usepackage{bm}

\usepackage[utf8]{inputenc}
\usepackage[T1]{fontenc}
\usepackage{mathptmx}
\usepackage{etoolbox}
\usepackage{sidecap}
\usepackage{fullpage}

\usepackage{grffile}
\graphicspath{{figures/}}

\usepackage{pgfplots}
\pgfplotsset{width=10cm,compat=1.9}
\usepackage{xcolor}
\definecolor{dark-violet}{RGB}{148,0,211}
\definecolor{sea-green}{RGB}{46, 139,  87}

\DeclareGraphicsExtensions{.pdf,.eps,.png,.jpg,.jpeg}
\usepackage{algorithm}
\usepackage{algpseudocode}

\newcommand{\bfalpha}{\boldsymbol \alpha}
\newcommand{\bfbeta}{\boldsymbol \beta}

\newcommand{\bff}{\boldsymbol f}

\newcommand{\Dcal}{\mathcal{D}}

\newcommand{\Tcal}{\mathcal{T}}

\newcommand{\Vcal}{\mathcal{V}}

\newcommand{\bfmu}{\boldsymbol \mu}

\newcommand{\bfe}{\boldsymbol e}

\newcommand{\bfp}{\boldsymbol p}
\newcommand{\bfq}{\boldsymbol q}

\newcommand{\bfs}{\boldsymbol s}

\newcommand{\bfA}{\boldsymbol A}

\newcommand{\bfC}{\boldsymbol C}

\newcommand{\bfF}{\boldsymbol F}

\newcommand{\bfP}{\boldsymbol P}
\newcommand{\bfR}{\boldsymbol R}

\newcommand{\bfV}{\boldsymbol V}

\newcommand{\bfQ}{\boldsymbol Q}

\newcommand{\bfS}{\boldsymbol S}

\newcommand{\nh}{N}
\newcommand{\nr}{n}
\newcommand{\nrp}{m}
\newcommand{\nrs}{m_s}

\newcommand{\win}{w}
\newcommand{\winit}{w_{\text{init}}}
\newcommand{\nz}{z}
\newcommand{\bbfS}{\breve{\bfS}}

\newcommand{\tbff}{\tilde{\bff}}

\newcommand{\tbfQ}{\tilde{\bfQ}}

\newcommand{\Qcal}{\mathcal{Q}}

\newcommand{\tbfq}{\tilde{\bfq}}
\newcommand{\hbfq}{\hat{\bfq}}

\newcommand{\kp}{l}

\newcommand{\hbff}{\hat{\bff}}

\newtheorem{remark}{Remark}

\newcommand{\onlinecite}[1]{\citenum{#1}}
\usepackage[superscript,biblabel]{cite}

\usepackage{lineno,hyperref}
\ifpdf
\hypersetup{
  pdftitle={Lookahead data-gathering strategies for online adaptive model reduction of transport-dominated problems}, 
  pdfauthor={Rodrigo Singh and Wayne Isaac Tan Uy and Benjamin Peherstorfer} 
}
\fi

\begin{document}

\title{Lookahead data-gathering strategies for online adaptive model reduction of transport-dominated problems}

\author{Rodrigo Singh\footnote{Courant Institute of Mathematical Sciences, New York University, 251 Mercer Street, New York, NY 10012, USA} \and Wayne Isaac Tan Uy\footnote{Analytics, Computing, and Complex Systems Laboratory, 123 Paseo de Roxas, Makati City, 1229, Philippines} \and Benjamin Peherstorfer\footnotemark[1]}

\date{\today}

\maketitle

\begin{abstract}
Online adaptive model reduction efficiently reduces numerical models of transport-dominated problems by updating reduced spaces over time, which leads to nonlinear approximations on latent manifolds that can achieve a faster error decay than classical linear model reduction methods that keep reduced spaces fixed. Critical for online adaptive model reduction is coupling the full and reduced model to judiciously gather data from the full model for adapting the reduced spaces so that accurate approximations of the evolving full-model solution fields can be maintained. In this work, we introduce lookahead data-gathering strategies that predict the next state of the full model for adapting reduced spaces towards dynamics that are likely to be seen in the immediate future. Numerical experiments demonstrate that the proposed lookahead strategies lead to accurate reduced models even for problems where previously introduced data-gathering strategies that look back in time fail to provide predictive models. The proposed lookahead strategies also improve the robustness and stability of online adaptive reduced models.
\end{abstract}

\section{Introduction}
Classical model reduction\cite{RozzaPateraSurvey,SIREVSurvey,SerkanInterpolatory} seeks linear approximations in reduced spaces and thus is inefficient for reducing numerical models that describe problems dominated by transport such as advecting flows and wave phenomena; see Ref.~\onlinecite{Notices} for a survey. This limitation of classical model reduction can be rigorously stated with the decay of the Kolmogorov $\nr$-width\cite{10.1093/imanum/dru066,MADAY2002289,OHLBERGER2013901,Greif19,Cohen2021} and thus is often referred to as the Kolmogorov barrier.
There is a range of numerical methods for nonlinear model reduction that circumvent the Kolmogorov barrier to construct efficient reduced models of transport-dominated problems.\cite{Notices} We broadly distinguish between methods that apply nonlinear transformations to recover linear low-dimensional structure and online adaptive methods that evolve nonlinear parametrizations in time. An early work that applies a nonlinear transformation to recover linear low-dimensional structure is Ref.~\onlinecite{ROWLEY20001}, which aims to account for symmetries. A range of other transformation have been developed, including transformations based on Wasserstein metrics,\cite{ehrlacher19} deep networks and deep autoencoders,\cite{LEE2020108973,KIM2022110841,https://doi.org/10.48550/arxiv.2203.00360} shifted proper orthogonal decomposition and its extensions,\cite{doi:10.1137/17M1140571,PAPAPICCO2022114687} quadratic manifolds,\cite{https://doi.org/10.48550/arxiv.2205.02304,BARNETT2022111348} and others.\cite{OHLBERGER2013901,TaddeiShock,https://doi.org/10.48550/arxiv.1911.06598,Cagniart2019,doi:10.1137/20M1316998} 

Online adaptive model reduction methods\cite{doi:10.1137/050639703,SAPSIS20092347,PhysRevE.89.022923,doi:10.1002/nme.4800,Peherstorfer15aDEIM,Amsallem2015,https://doi.org/10.1002/nme.4770,P18AADEIM,RAMEZANIAN2021113882} evolve nonlinear parametrizations over time. We build on ADEIM\cite{Peherstorfer15aDEIM,P18AADEIM,ZPW17SIMAXManifold,CKMP19ADEIMQuasiOptimalPoints,UWHP22Flame} that adapts basis functions with low-rank updates. There is a close connection to dynamic low-rank approximations and related techniques that evolve time-dependent nonlinear parametrizations\cite{doi:10.1137/050639703,SAPSIS20092347,doi:10.1137/16M1088958,Pagliantini2021,Peherstorfer2016,doi:10.1137/140967787,Musharbash2020,refId0Hesthaven,hesthaven_pagliantini_rozza_2022,https://doi.org/10.48550/arxiv.2203.01360,wen2023coupling} as well as to the Dirac-Frenkel variational principle.\cite{dirac_1930,Frenkel1934,Lubich2008,lasser_lubich_2020}. A distinguishing feature of ADEIM is that it remains efficient even when reducing numerical models with nonlinear state dynamics because ADEIM uses empirical interpolation\cite{Everson1995,barrault_empirical_2004,AstWWB08,deim2010,QDEIM,PDG18ODEIM} to approximate the basis updates from only sparse sketches of approximate full-model states. In Ref.~\onlinecite{CKMP19ADEIMQuasiOptimalPoints}, a sampling scheme is introduced that quasi-optimally selects the state components to sketch for adapting reduced spaces. In Ref.~\onlinecite{UWHP22Flame}, online adaptive model reduction with ADEIM has been applied to derive reduced models of benchmark problems motivated by chemically reacting flows.\cite{10.1007/3-540-31801-1_89,chris_wentland_2021_5517532,10.1007/3-540-31801-1_89,HUANGMPLSVT} The authors of Ref.~\onlinecite{huang2023predictive} introduce basis updates that are motivated by Ref.~\onlinecite{ZPW17SIMAXManifold} and preserve properties when updating the basis such as orthogonality. Ref.~\onlinecite{zucatti2023adaptive} uses ADEIM as a building block to derive an online adaptive method with a sub-iteration scheme that solves the full model on a sparse sketch to adapt the reduced space. The focus of Ref.~\onlinecite{zucatti2023adaptive} is on problems with shocks.

Critical for online adaptive model reduction methods such as ADEIM is the coupling of the full and reduced model to judiciously gather data for adapting the reduced space so that the reduced model can maintain an accurate approximation of the evolving full-model solution field over time. We introduce lookahead strategies that derive reduced predictor models that are cheap to simulate and that determine at which state to sketch the full model. The predictor models are set up such that the sketch of the approximate full-model state corresponds to dynamics that will be likely visited in the immediate future so that a sketch provides informative data for adapting the reduced spaces. This is in contrast to previous strategies\cite{P18AADEIM} that gather sketches from the full model that lag behind in time; which is also observed in Ref.~\onlinecite{zucatti2023adaptive}. We demonstrate on numerical examples that the lookahead strategies lead to more accurate approximations than previous data-gathering strategies that look back in time. The lookahead strategies also improve the robustness of ADEIM reduced models, which is essential for long-time integration.

\section{Preliminaries}\label{sec:Prelim}
We briefly discuss online adaptive model reduction with ADEIM in Sections~\ref{sec:Prelim:FOM}--\ref{sec:Prelim:UpdateSamplingPoints} and then formulate the problem in Section~\ref{sec:Prelim:Problem}.

\subsection{Numerical models}\label{sec:Prelim:FOM}
A widely used approach to numerically solve systems of partial differential equations (PDEs) is to discretize them in space with, e.g., finite-volume, finite-element, and finite-difference methods to obtain dynamical-system models (systems of ordinary differential equations), which are then numerically integrated in time.
We denote such a dynamical-system model as  
    \begin{equation}\label{eq:contFOM}
        \frac{\mathrm d}{\mathrm dt} \bfq(t;\bfmu) = \bar{\bff}(\bfq(t;\bfmu); \bfmu)\,,\qquad t \in (0,T],
    \end{equation}
    where $\bfmu \in \Dcal$ is a physical parameter, $\bfq(t;\bfmu) \in \mathbb{R}^{\nh}$ is the state at time $t$ and $\bar{\bff}: \mathbb{R}^{\nh} \times \Dcal \to \mathbb{R}^{\nh}$ is a nonlinear vector-valued function that encodes the dynamics of the problem. 

Applying an implicit time discretization to \eqref{eq:contFOM} with time-step size $\delta t > 0$ leads to the time-discrete model 
    \begin{equation}\label{eq:discFOM}
        \bfq_{k-1}(\bfmu) = \bff(\bfq_k(\bfmu); \bfmu)\,,\qquad k = 1, \dots, K\,,
    \end{equation}
    where $\bff: \mathbb{R}^{\nh} \times \Dcal \to \mathbb{R}^{\nh}$ denotes the time-discrete dynamics and $\bfq_k(\bfmu)$ is the approximation of $\bfq(t_k; \bfmu)$ at time step $t_k = \delta t k$ for $ k = 0, \dots, K$. 
In the following, we are interested in situations where fine discretizations in space are necessary to resolve the phenomena of interest, which results in high dimensions $\nh$ of states $\bfq(t; \bfmu)$ of model \eqref{eq:contFOM} and high computational costs of time stepping the time-discrete model \eqref{eq:discFOM}.

\subsection{Static projection-based reduced models}\label{sec:Prelim:StaticROM}
Model reduction \cite{RozzaPateraSurvey,SIREVSurvey,SerkanInterpolatory} is typically split into a one-time high-cost offline phase where data are collected from a full model to build a reduced model, and an online phase where the reduced model is used to compute reduced state trajectories that approximate the full-model state trajectories at new parameters and new initial conditions.
We start by describing the offline phase of classical model reduction: Let $\bfmu_1, \dots, \bfmu_M \in \Dcal$ be training parameters and $\bfq_0(\bfmu_1), \dots, \bfq_0(\bfmu_M)$ be initial conditions for simulating the full model \eqref{eq:discFOM}. 
The snapshot matrix is $\bfQ = [\bfQ(\bfmu_1), \dots, \bfQ(\bfmu_M)] \in \mathbb{R}^{\nh \times M(K+1)}$, where the training trajectories
    \begin{equation}
    \bfQ(\bfmu_i) = [\bfq_0(\bfmu_i), \dots, \bfq_K(\bfmu_i)]\,,\qquad i = 1, \dots, M\,,
    \end{equation}
    are obtained by time stepping the full model \eqref{eq:discFOM}.
Applying proper orthogonal decomposition (POD) \cite{RozzaPateraSurvey,SIREVSurvey} to the snapshot matrix $\bfQ$ yields an orthonormal basis matrix $\bfV \in \mathbb{R}^{\nh \times \nr}$ whose columns span a subspace $\Vcal$ of $\mathbb{R}^{\nh}$ of dimension $\nr < \nh$.
Projecting the model operators of the full model \eqref{eq:discFOM} onto the subspace $\Vcal$ leads to the Galerkin reduced model  
\begin{equation}\label{eq:Prelim:GalerkinROM}
    \tilde{\bfq}_{k-1}(\bfmu) = \bfV^T\bff(\bfV\tilde{\bfq}_k(\bfmu); \bfmu)\,,\qquad k = 1, \dots, K\,,
    \end{equation}
    with the reduced state $\tbfq_k(\bfmu) \in \mathbb{R}^{\nr}$ at time step $k$ of dimension $\nr$. The initial condition is projected as $\tbfq_0(\bfmu) = \bfV^T \bfq_0(\bfmu) \in \tilde{\Qcal}_0 \subseteq \mathbb{R}^{\nr}$. The lifted reduced state $\bfV \tbfq_k(\bfmu)$ is an approximation of $\bfq_k(\bfmu)$.
Except in limited cases, e.g., if $\bff$ is polynomial in the state $\bfq$, solving the Galerkin reduced model \eqref{eq:Prelim:GalerkinROM} requires evaluating $\bff$ at all components such that the computational costs scale with the dimension $\nh$. 
To circumvent the scaling of the costs with the dimension $\nh$ of the full-model states, empirical interpolation \cite{barrault_empirical_2004,grepl_efficient_2007,deim2010,QDEIM} constructs an approximation $\tbff: \mathbb{R}^{\nr} \times \Dcal \to \mathbb{R}^{\nr}$ of $\bff$ of the form
    $\tbff(\tbfq; \bfmu) = (\bfP^T\bfV)^{\dagger}\bfP^T\bff(\bfV\tilde{\bfq}; \bfmu)$
    where $(\bfP^T\bfV)^{\dagger}$ denotes the Moore--Penrose pseudo-inverse of $\bfP^T\bfV$. 
The matrix $\bfP \in \mathbb{R}^{\nh \times \nrp}, \nrp \ge \nr$ selects $\nrp$ component functions of $\bff$ that are evaluated:
    It is defined as $\bfP = [\bfe_{p_1}, \dots, \bfe_{p_\nrp}]\in \{0,1\}^{\nh \times \nrp}$, with pairwise disjoint indices $p_1,\dots,p_m \in \{1,\dots,N\}$. 
     For $i=1,\dots,\nrp$, the vector $\bfe_{p_i} \in \{0,1\}^{\nh}$ is the canonical unit vector with 1 at the component indexed by $p_i$ and zero at all other components. We define the vector $\bfp = [p_1, \dots, p_{\nrp}]^T$.
Multiplication with the matrix $\bfP^T$ selects $\nrp$ components from an $\nh$-dimensional vector so that $\bfP^T\bff(\bfV \tbfq; \bfmu)$ requires evaluating only the $\nrp$ component functions of $\bff$ corresponding to the indices $p_1, \dots, p_{\nrp}$.  
Approaches to compute the indices $p_1,\dots,p_{\nrp}$ include greedy methods, \cite{barrault_empirical_2004,deim2010} QDEIM, \cite{QDEIM} and oversampling algorithms, \cite{PDG18ODEIM,doi:10.2514/6.2021-1371}.
Combining the Galerkin reduced model with empirical interpolation yields the Galerkin-EIM reduced model 
    \begin{equation}\label{eq:DEIMROM}
        \tilde{\bfq}_{k-1}(\bfmu) = \tilde{\bff}(\tilde{\bfq}_k(\bfmu); \bfmu)\,,\qquad k = 1, \dots, K\,.
    \end{equation}
Depending on the structure\cite{barrault_empirical_2004,deim2010} of $\bff$, the costs of time stepping \eqref{eq:DEIMROM} scale with the reduced dimension $\nr$ and the number of indices $\nrp$ but are independent of the dimension $\nh$ of the states of the full model.
Once a reduced model has been constructed in the offline phase, it is used in the online phase to compute reduced states $\tbfq_1(\bfmu), \dots, \tbfq_K(\bfmu)$ with a new initial condition $\tbfq_0(\bfmu)$ and parameter $\bfmu \in \Dcal \setminus \{\mu_1, \dots, \mu_M\}$.
Speedups compared to simulating the full model \eqref{eq:discFOM} can be achieved if the one-time high-cost offline phase is compensated by evaluating the reduced model at many parameters online.

\subsection{Online adaptive model reduction with ADEIM}\label{sec:Prelim:ADEIM}
Online adaptive reduced models based on ADEIM\cite{Peherstorfer15aDEIM,P18AADEIM,CKMP19ADEIMQuasiOptimalPoints,ZPW17SIMAXManifold,UWHP22Flame} update the basis matrix $\bfV$ and thus the subspace $\Vcal$ during the online phase. The online adaptation breaks\cite{Peherstorfer15aDEIM} with the classical offline/online decomposition of first constructing and then using the reduced model. 

\begin{remark}
Online adaptive reduced models based on ADEIM\cite{Peherstorfer15aDEIM} are training-free in the sense that there is no one-time, high-cost offline phase. In particular, an ADEIM reduced model is adapted towards the current parameter $\bfmu$ of interest and thus does \emph{not} require extensively sampling the parameter domain $\Dcal$, which is challenging in high dimensions.
\end{remark}

Consider the Galerkin-EIM model defined in \eqref{eq:DEIMROM} and the approximation $\tbff$. The approximation $\tbff$ depends on the basis matrix $\bfV$ and on the selection matrix $\bfP$. 
Let the basis matrix $\bfV_k$ and the selection matrix $\bfP_k$ depend on the time step $k = 0, \dots, K$ so that the approximation $\tbff_k$ also depends on time,
\begin{equation}\label{eq:Prelim:FTilde}
\tbff_k(\tbfq_k(\bfmu); \bfmu) = (\bfP_k^T\bfV_k)^{\dagger}\bfP_k^T\bff(\bfV_k\tbfq_k(\bfmu); \bfmu)\,,
\end{equation}
The corresponding online adaptive reduced model is
\begin{equation}\label{eq:Prelim:ADEIMModel}
\tbfq_{k-1}(\bfmu) = \tbff_k(\tbfq_k(\bfmu); \bfmu)\,,\qquad k = 1, \dots, K\,.
\end{equation}
We build on ADEIM to adapt the basis matrix at time step $k$ from $\bfV_k$ and space $\Vcal_k$ to $\bfV_{k+1}$ and $\Vcal_{k + 1}$. The adaptation is $\bfV_{k+1} = \bfV_k + \bfalpha_k\bfbeta_k^T$ where $\bfalpha_k \in \mathbb{R}^{\nh}$ and $\bfbeta_k \in \mathbb{R}^{\nr}$ provide a rank-1 update.
To compute the update $\bfalpha_k\bfbeta_k^T$, ADEIM defines a data window of $w \in \mathbb{N}$ data samples,
\begin{equation}\label{eq:Prelim:UpdateWindow}
\bfF_k = [\hbfq_{k - w + 1}(\bfmu), \dots, \hbfq_k(\bfmu)] \in \mathbb{R}^{\nh \times \win}\,.
\end{equation}
The update $\bfalpha_k\bfbeta_k^T$ is the solution to the optimization problem 
\begin{equation}
\min_{\bfalpha_k \in \mathbb{R}^{\nh},\, \bfbeta_k \in \mathbb{R}^{\nr}}\, \left\|(\bfV_k + \bfalpha_k\bfbeta_k^T)\bfC_k - \bfF_k\right\|_F^2\,,\label{eq:ADEIMUpdate}
\end{equation}
where $\bfC_k = \bfV_k^T\bfF_k$ is the coefficient matrix.
The update $\bfalpha_k\bfbeta_k^T$ computed with \eqref{eq:ADEIMUpdate} adapts the space with the aim that $\Vcal_{k + 1}$ can approximate well the data points in $\bfF_k$ in the sense that the projection error $\|\bfV_{k + 1}\bfV_{k + 1}^{\dag}\bfF_k - \bfF_k\|_F$ is low. As a side remark, we state that the optimization problem \eqref{eq:ADEIMUpdate} is different from the original ADEIM problem\cite{Peherstorfer15aDEIM,P18AADEIM} because the norm of the residual in \eqref{eq:ADEIMUpdate} is taken over all components rather than over a sparse subset of the set of all components.

\subsection{Data collection for adaptation}\label{sec:Prelim:Lookback}
The choice of the data samples in $\bfF_k$ is critical because the space is adapted such that it approximates well the data samples. 
It is  proposed in Ref.~\onlinecite{P18AADEIM} to use the following insight for deriving data samples for the adaptation: if the lifted reduced state $\bfV_k\tbfq_k(\bfmu)$ of the reduced model \eqref{eq:Prelim:ADEIMModel} at time step $k$ approximates well the state $\bfq_k(\bfmu)$ of the full model \eqref{eq:discFOM}, then evaluating the full-model right-hand side function $\bff$ at the lifted reduced state will give an approximation $\bff(\bfV_k\tbfq_k(\bfmu))$ of the full-model state $\bfq_{k - 1}(\bfmu)$ at time step $k - 1$. This motivates using the approximation $\bff(\bfV_k\tbfq_k(\bfmu))$ of the full-model state $\bfq_k(\bfmu)$ to inform the data sample $\hbfq_k(\bfmu)$.

Evaluating $\bff$ is typically computationally cheaper than taking a time step with the full model \eqref{eq:discFOM}, which requires a (nonlinear) solve. However, the costs of evaluating the full-model right-hand side function $\bff$ still scales with the dimension $\nh$ of the full-model states. To avoid this, it is proposed\cite{P18AADEIM} to evaluate $\bff$ only at a sketch of $\nrs \ll \nh$ components based on a sampling matrix $\bfS_k \in \{0, 1\}^{\nh \times \nrs}$. 
The sampling matrix $\bfS_k$ is defined analogously to the selection matrix $\bfP_k$ by defining the vector $\bfs_k = [s_1^{(k)}, \dots, s_{\nrs}^{(k)}]^T \in \{1, \dots, \nh\}^{\nrs}$ of pairwise distinct indices of components functions of $\bff$ that give rise to a matrix 
$\bfS_k = [\bfe_{s_1^{(k)}}, \dots, \bfe_{s_{\nrs}^{(k)}}]$
via the canonical unit vectors of dimension $\nh$.

With the sampling matrix $\bfS_k$, at time step $k$, the full-model right-hand side function $\bff$ is evaluated only at the component functions corresponding to the sampling points $s_1^{(k)}, \dots, s_{\nrs}^{(k)}$ 
and all other components are approximated with empirical interpolation based on the current subspace $\Vcal_k$. This leads to the data sample $\hbfq_k(\bfmu)$ at time $k$ defined as
\begin{equation}\label{eq:Prelim:LookBack}
\begin{aligned}
\bfS_k^T\hbfq_k(\bfmu) = &  \bfS_k^T\bff(\bfV_k\tbfq_k(\bfmu))\,,\qquad \\
\bbfS_k^T\hbfq_k(\bfmu) = &\bbfS_k^T\bfV_k(\bfS_k^T\bfV_k)^{\dagger}\bfS_k^T\bff(\bfV_k\tbfq_k(\bfmu))\,,
\end{aligned}
\end{equation}
where $\bbfS_k \in \{0, 1\}^{\nh \times (\nh - \nrs)}$ is the complementary sampling point matrix to $\bfS_k$ corresponding to the sampling points in $\{1, \dots, \nh\}\setminus \{s_1^{(k)}, \dots, s_{\nrs}^{(k)}\}$. We note that several works\cite{P18AADEIM,UWHP22Flame} use different combinations of the points corresponding to $\bfP_k$ and $\bfS_k$ to compute the approximation \eqref{eq:Prelim:LookBack}.

Once the basis matrix $\bfV_k$ is adapted to $\bfV_{k + 1}$, then the selection matrix $\bfP_{k+1}$ is adapted by applying either QDEIM \cite{QDEIM} or ODEIM\cite{PDG18ODEIM} onto the adapted basis matrix $\bfV_{k + 1}$; the costs of which scale linearly with the dimension $\nh$.
The adapted basis matrix $\bfV_{k + 1}$ and selection matrix $\bfP_{k + 1}$ (in vector notation $\bfp_{k + 1}$) define $\tbff_{k + 1}$ as in \eqref{eq:Prelim:FTilde} for solving the online adaptive reduced model \eqref{eq:Prelim:ADEIMModel} at time step $k + 1$.

\subsection{Updating the sampling points in ADEIM}\label{sec:Prelim:UpdateSamplingPoints}
As shown in Ref.~\onlinecite{P18AADEIM}, the sampling points corresponding to the sampling matrix $\bfS_k$ have to be updated for approximating well transport-dominated problems with ADEIM. 
We will update the sampling matrix $\bfS_k$ every $z$-th time step, where $z \in \mathbb{N}$. 
To update $\bfS_k$ to $\bfS_{k + 1}$, the $\nh \times \win$ residual matrix is computed as
$\bfR_k = \bfF_k - \bfV_k(\bfP_k^T\bfV_k)^{\dagger}\bfP_k^T\bfF_k$
where $\bfF_k$ is the data window at time step $k$, the matrix $\bfP_k$ encodes the indices for the empirical interpolation, and $\bfV_k$ is the basis matrix.
The approach proposed in Ref.~\onlinecite{P18AADEIM} to update the sampling points is to compute the Euclidean norms of the rows of the matrix $\bfR_k$, which we denote as $r_k^{(1)}, \dots, r_k^{(\nh)}$, sorting them in descending order $r_k^{(i_1)} \geq r_k^{(i_2)} \geq \cdots \geq r_k^{(i_{\nh})} \geq 0$ and then setting the first $\nrs$ indices $i_1, \dots, i_{\nrs}$ as new sampling points $s_1^{(k+1)} = i_1, \dots, s_{\nrs}^{(k+1)} = i_{\nrs}$, with the corresponding sampling matrix $\bfS_{k + 1}$.
This process leads to quasi-optimal sampling points.\cite{CKMP19ADEIMQuasiOptimalPoints}

\subsection{Problem formulation}\label{sec:Prelim:Problem}
The data window $\bfF_k$ defined in \eqref{eq:Prelim:UpdateWindow} drives the adaptation of the reduced space $\Vcal_k$. 
Based on the procedure described in Section~\ref{sec:Prelim:ADEIM}, at time step $k$, the reduced state $\tbfq_k(\bfmu)$ is used to evaluate the full-model right-hand side function $\bff$ to obtain the data sample $\hbfq_k(\bfmu)$ defined in \eqref{eq:Prelim:LookBack}, which approximates $\bfq_{k - 1}(\bfmu)$. Thus, the data sample $\hbfq_k(\bfmu)$ approximates the full-model state $\bfq_{k - 1}(\bfmu)$ at time step $k - 1$, which means that the data window $\bfF_k$ is filled with information from the previous time step $k - 1$ to update the space at the next time step $k$. This means that the approach\cite{P18AADEIM} of Section~\ref{sec:Prelim:ADEIM} looks back in time to fill the data window $\bfF_k$ and thus the adaptation of the space lags behind.

\section{Lookahead strategies for online adaptive model reduction}\label{sec:Ahead}
We introduce lookahead strategies for generating data samples from the full model that correspond to dynamics that will be likely seen in the immediate future for online adaptive model reduction. By looking ahead, the data samples in the data window are informative about what dynamics will be seen in the future for adapting reduced spaces, rather than what has been seen in the past. The looking ahead also helps to improve stability of adaptive reduced models in our numerical experiments.  

For ease of exposition, we drop the dependence on the parameter $\bfmu$ in the notation in most of this section.

\subsection{Predictor model for looking ahead in time}
Consider a time-step size $\delta \tau = \delta t/C_{\tau}$, where $\delta t$ is the time-step size of the full model \eqref{eq:discFOM} and $C_{\tau} \in \mathbb{N}$ is an integer factor.
Define a counter variable $\kp = 1, 2, \dots, $ with time steps $0 = \tau_0 < \tau_1 < \tau_2 < \dots < \tau_{\kp} < \dots $ so that $\tau_{k C_{\tau}} = t_k$ for $k = 1, \dots, K$. 
We now define a predictor model that is obtained by discretizing the time-continuous full model \eqref{eq:contFOM} with an explicit time integration scheme and time-step size $\delta \tau$,
\begin{equation}\label{eq:Ahead:PredictorModel}
\bfq_{\kp}^{(P)} = \bff^{(P)}(\bfq_{\kp - 1}^{(P)})\,,\qquad \kp = 1, 2, 3, \dots\,,
\end{equation}
where $\bfq_{\kp}^{(P)}$ is the state at time step $\kp$ and $\bff^{(P)}$ is the right-hand side function obtained after time discretization with $\delta \tau$. 
Taking a time step with the predictor model \eqref{eq:Ahead:PredictorModel} incurs costs that scale with the costs of evaluating the right-hand side function $\bff^{(P)}$ because an explicit time-integration scheme is used in \eqref{eq:Ahead:PredictorModel}, which typically incurs lower costs than solving the nonlinear system of equations for taking a time step with the full model \eqref{eq:discFOM} that is based on an implicit time integration scheme. 
However, solving the predictor model \eqref{eq:Ahead:PredictorModel} still incurs costs that scale at least linearly in the dimension $\nh$ of the full-model states. Furthermore, explicit schemes can require small time-step sizes $\delta \tau$, which can further increase the computational costs.

\subsection{Reducing the costs of simulating the predictor model}
We introduce a reduced predictor model that uses empirical interpolation to rapidly approximate the right-hand side function $\bff^{(P)}$ of the predictor model \eqref{eq:Ahead:PredictorModel} with the current basis matrix $\bfV_k$ and sampling matrix $\bfS_k$ local in time\,,
\begin{equation}\label{eq:Ahead:ReducedPredictorModel}
\hbfq_{\kp}^{(P)} = \hbff^{(P)}_{k}(\hbfq_{\kp - 1}^{(P)})\,,\qquad \kp = 1, 2, 3, \dots\,.
\end{equation}
The approximation $\hbff^{(P)}_{k}$ is 
\begin{equation}\label{eq:PredictorEIM}
\begin{aligned}
\bfS_k^T\hbff^{(P)}_{k}(\hbfq_{\kp - 1}^{(P)}) = & \bfS_k^T\bff^{(P)}(\hbfq_{\kp - 1}^{(P)})\,,\\
\bbfS_k^T\hbff^{(P)}_{k}(\hbfq_{\kp - 1}^{(P)}) = & \bbfS_k^T\bfV_k(\bfS_k^{T}\bfV_k)^{\dagger}\bfS_k^T\bff^{(P)}(\hbfq_{\kp - 1}^{(P)})\,,
\end{aligned}
\end{equation}
which sets the component of $\hbff_k^{(P)}$ selected by $\bfS_k$ to the values of $\bff^{(P)}$ and approximates all other components given by the complementary sampling matrix $\bbfS_k$ with empirical interpolation. 
The approximation \eqref{eq:PredictorEIM} is fixed (non-adaptive) with respect to the time steps $\kp = 1, 2, 3, \dots$ of the predictor model but changes with the the time step $k$ of the adaptive reduced model.

To set the costs of time stepping the reduced predictor model \eqref{eq:PredictorEIM} in context, recall that taking a time step with the full model \eqref{eq:discFOM} typically incurs the costs of solving a nonlinear system of $\nh$ equations. 
Taking a time step with the predictor model \eqref{eq:Ahead:PredictorModel} incurs the costs of evaluating the full-model right-hand side function at $\nh$ components, which is in contrast to the full model \eqref{eq:discFOM} because the predictor model uses an explicit time-stepping scheme so the costs per time step are lower. 
The reduced predictor model \eqref{eq:Ahead:ReducedPredictorModel} evaluates the full-model right-hand side function $\bff^{(P)}$ at only $\nrs$ components instead of at all $\nh$ components and thus we expected the costs of a time step to be lower than with predictor model \eqref{eq:Ahead:PredictorModel} and the full model \eqref{eq:discFOM}.

\begin{remark}\label{remark:UseVInPredictoROM}
Instead of just approximating the right-hand side function $\bff^{(P)}$ with empirical interpolation as in the reduced predictor model \eqref{eq:PredictorEIM}, we could have also approximated the state $\bfq_{\kp}^{(P)}$ of the predictor model \eqref{eq:Ahead:PredictorModel} in the reduced space $\bfV_k$. Such an additional approximation would lead to a reduced predictor model with costs per time step that scale independently of the dimension $\nh$, if $\bff^{(P)}$ has sufficient structure for empirical interpolation.\cite{barrault_empirical_2004,deim2010}
\end{remark}

\subsection{Lookeahead strategy for filling data window}
We now use the reduced predictor model \eqref{eq:Ahead:ReducedPredictorModel} to compute a new data sample $\hbfq_k$ for the data window $\bfF_k$ at time step $k$. 
Because the costs of taking a time step with the reduced predictor model \eqref{eq:Ahead:ReducedPredictorModel} are low compared to taking a time step with the full model, we propose to used the lifted reduced state $\bfV_k\tbfq_k$ at time step $k$ as initial condition for the reduced predictor model \eqref{eq:Ahead:ReducedPredictorModel} and simulate it for $C_{\tau}$-many time steps to get an approximation of the full model state at time step $k + 1$, which is then used in the data window $\bfF_k$ to inform a new data sample $\hbfq_k$. 
By using the reduced predictor model and integrating it in time for a few time steps, we obtain a prediction of the state at the subsequent time step $k + 1$ and thus this strategy looks ahead in time.

Given the reduced state $\tbfq_k$ at time step $k$, the reduced predictor model \eqref{eq:Ahead:ReducedPredictorModel} is integrated in time for $C_{\tau}$-many time steps starting from the initial condition $\hbfq_0^{(P)} = \bfV_k\tbfq_k$ to obtain $\hbfq^{(P)}_{C_\tau}$. The state $\hbfq^{(P)}_{C_\tau}$ then serves as the data sample $\hbfq_k$ at time step $k$ for filling the data window $\bfF_k$.

Notice that now $\hbfq_k$ is an approximation of the full-model state $\bfq_{k + 1}$ at time step $k + 1$.
At the same time, notice that taking a time-step with the full model \eqref{eq:discFOM} requires numerically solving a system of nonlinear equations, which is avoided and thus computing the data sample $\hbfq_k$ for the data window has costs that scale independently of the costs of taking a time step with the full model \eqref{eq:discFOM}.

\begin{remark}
Building on Remark~\ref{remark:UseVInPredictoROM}, we note that if the basis $\bfV_k$ is used to approximate also the state of the reduced predictor model as $\tbfq_{\kp}^{(P)} \in \mathbb{R}^{\nr}$, then a data sample $\hbfq_k$ can be obtained by evaluating the right-hand side function $\bff^{(P)}$ of the predictor model \eqref{eq:Ahead:PredictorModel} at the lifted reduced state $\bfV_k\tbfq_{C_{\tau}}^{(P)}$ after $C_{\tau}$ time steps as
\begin{equation}
\begin{aligned}
\bfS_k^T\hbfq_k = & \bfS_k^T\bff^{(P)}(\bfV_k\tbfq^{(P)}_{C_{\tau}})\,,\\
\bbfS_k^T\hbfq_k = &\bbfS_k^T\bfV_k(\bfS_k^T\bfV_k)^{\dagger}\bfS_k^T\bff^{(P)}(\bfV_k\tbfq^{(P)}_{C_{\tau}})\,.
\end{aligned}
\end{equation}
Note that the lifted reduced predictor state $\bfV_k\tbfq_{C_{\tau}}^{(P)}$ is in the reduced space $\Vcal_k$ and thus carries no new information to adapt $\Vcal_k$ if used directly as data sample $\hbfq_k$. 
\end{remark}

\algrenewcommand\algorithmicindent{0.5em}
\begin{algorithm}[t]
\begin{algorithmic}[1]
\Procedure{ADEIM}{$\bfq_0, \bfmu, \nr, \winit, \win, \nrs, \nz, C_{\tau}$}
\State Iterate \eqref{eq:discFOM} for $\winit$ time steps and store states in $\bfQ$\label{alg:ABAS:SolveFOM}
\State Set $k = \winit+1$\label{alg:ABAS:StartROMInit}
\State Compute $\nr$-dimensional POD basis $\bfV_k$ of $\bfQ$\label{alg:ABAS:PODBasisConstruction}
\State Compute QDEIM interpolation points $\bfp_k$ for basis $\bfV_k$
\State Initialize data window $\bfF = \bfQ[:, k-\win+1:k-1]$
\State Initialize state $\tilde{\bfq}_{k - 1} = \bfV_k^T\bfQ[:, k - 1]$\label{alg:ABAS:EndROMInit}
\For{$k = \winit + 1, \dots, K$}\label{alg:ABAS:ROMLoop}
\State Solve \eqref{eq:Prelim:ADEIMModel} for $\tilde{\bfq}_k$, using basis matrix $\bfV_k$ and points $\bfp_k$\label{alg:ABAS:SolveROM}
\State Store lifted reduced state $\bfQ[:, k] = \bfV_k\tilde{\bfq}_k$
\If{$\operatorname{mod}(k, \nz) == 0 || k == \winit + 1$}\label{alg:ABAS:IfSampling}
\State Take $C_{\tau}$ time steps with full predictor \eqref{eq:Ahead:PredictorModel} with $\hbfq_0^{(P)} = \bfV_k\tbfq_k$
\State Store data sample $\bfF[:, k] = \hbfq_{C_{\tau}}^{(P)}$
\State $\bfR_k = \bfF[:, k - \win + 1:k] - \bfV_k(\bfV_k[\bfp_k, :])^{\dagger}\bfF[\bfp_k, k - \win + 1:k]$
\State Update sampling points $\bfs_k$ as in Section~\ref{sec:Prelim:UpdateSamplingPoints}
\Else
\State Take $C_{\tau}$ time steps with reduced predictor \eqref{eq:Ahead:ReducedPredictorModel}, $\hbfq_0^{(P)} = \bfV_k\tbfq_k$
\State Store data sample $\bfF[:, k] = \hbfq_{C_{\tau}}^{(P)}$
\EndIf
\State Solve \eqref{eq:ADEIMUpdate} for $\bfalpha_k,\bfbeta_k$ with $\bfF_k = \bfF[:, k - \win + 1:k]$ and $\bfV_k$\label{alg:ABAS:CompBasisUpdate}
\State Adapt basis $\bfV_{k + 1} = \bfV_k + \bfalpha_k\bfbeta_k$ and orthogonalize $\bfV_{k + 1}$\label{alg:ABAS:ApplyBasisUpdate}
\State Compute points $\bfp_{k + 1}$ by applying QDEIM to $\bfV_{k + 1}$ \label{alg:ABAS:AdaptP}
\EndFor\\
\Return Return trajectory $\bfQ$
\EndProcedure
\end{algorithmic}
\caption{ADEIM algorithm with lookahead}\label{alg:ADEIMLookahead}
\end{algorithm}

\begin{figure*}
\begin{tabular}{ccc}
\resizebox{0.30\textwidth}{!}{\huge\input{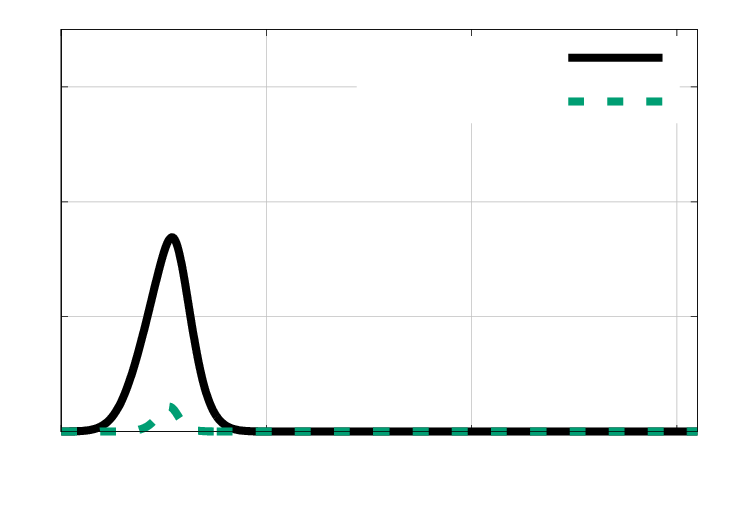}} & \resizebox{0.30\textwidth}{!}{\huge\input{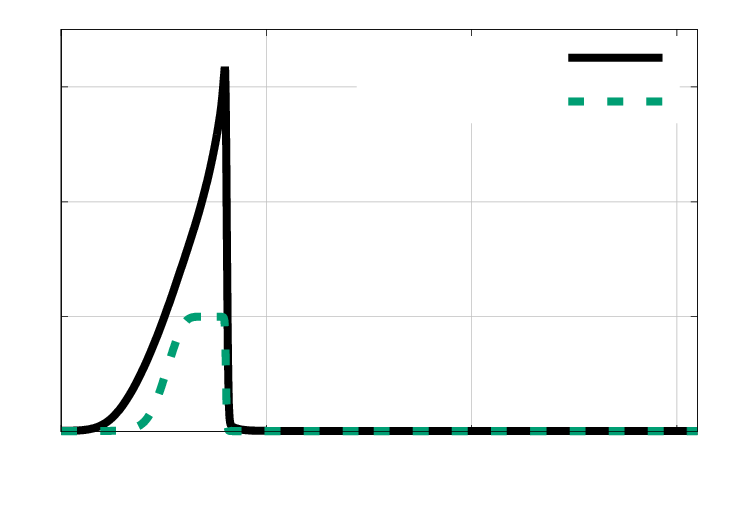}} & \resizebox{0.30\textwidth}{!}{\huge\input{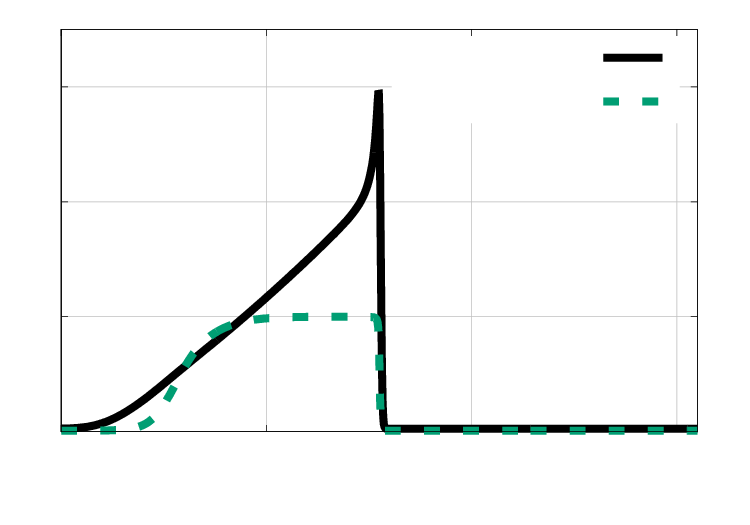}}\\
(a) time $t = 0.05$ & (b) time $t=0.25$ & (c) time $t=1$\\
\resizebox{0.30\textwidth}{!}{\huge\input{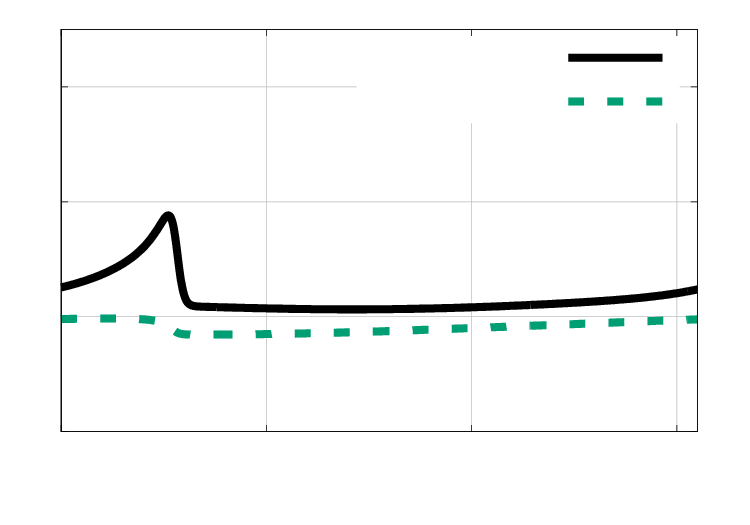}} & \resizebox{0.30\textwidth}{!}{\huge\input{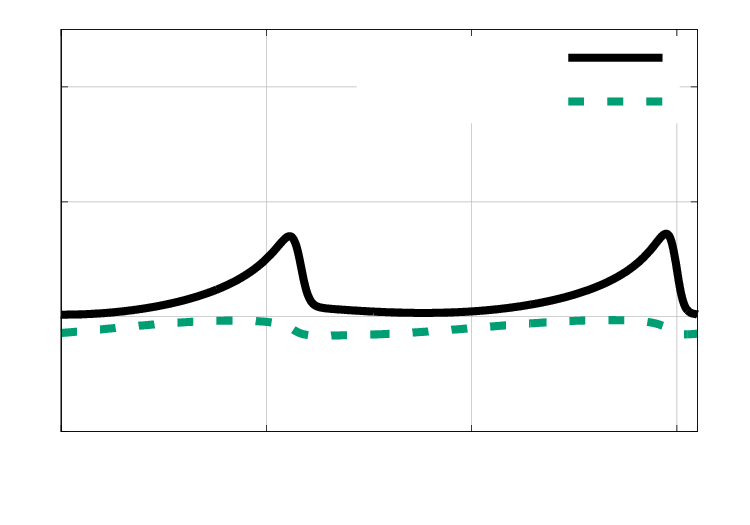}} & \includegraphics[width=0.33\textwidth]{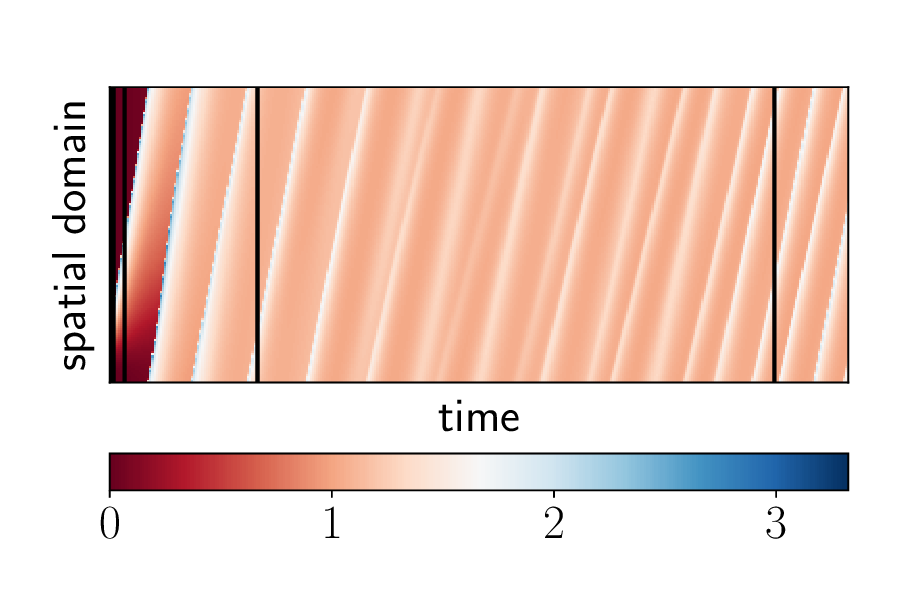}\\
(d) time $t=10$ & (e) time $t=45$ & (f) space-time plot of intensity $\eta$\\
\end{tabular}
\caption{Rotating detonation waves: Plots of the intensive property $\eta$ of the working fluid and the combustion progress $\lambda$ at times $t = 0.05, 0.25, 1, 10, 45$. The plots show the evolution of the smooth single-pulse initial condition into a solution close to a shock and then a wave that circulates and eventually spawns a second wave. Model based on Ref.~\onlinecite{PhysRevE.101.013106}.} 
\label{fig:RDE_FOM2DStills}
\end{figure*}

\subsection{Computational procedure}
We summarize the computational procedure in Algorithm~\ref{alg:ADEIMLookahead}. In the algorithm, we use canonical slicing notation, where $\bfA[i, :]$ and $\bfA[:, j]$ denotes the $i$-th row and the $j$-th column, respectively, of a matrix $\bfA$. If $\bfp \in \{1, \dots, \nh\}^{m}$ is an $m$-dimensional vector and the matrix $\bfA$ has $\nh$ many rows, then $\bfA[\bfp, :]$ selects the $m$ rows corresponding to the indices in the vector $\bfp$.

The inputs to the algorithm are the initial condition $\bfq_0 \in \mathbb{R}^{\nh}$ and the parameter $\bfmu \in \Dcal$ at which an approximation of the full-model trajectory should be computed. Other inputs are the dimension $\nr$ of the reduced model, the initial window size $\winit \in \mathbb{N}$, the data window size $\win$, the number of sampling points $\nrs$, the frequency $\nz$ with which to update the sampling points, and the factor $C_{\tau}$ of the time-step size $\delta \tau = \delta t / C_{\tau}$ of the reduced predictor model \eqref{eq:Ahead:ReducedPredictorModel}. 
In the first \ref{alg:ABAS:EndROMInit} lines, the reduced model is initialized by time-stepping the full model for $\winit \ll K$ time steps to compute snapshots for constructing the POD basis $\bfV_{\winit}$ and the points for empirical interpolation $\bfp_{\winit}$.
Starting with time step $k = \winit + 1$, the ADEIM reduced model is used for time stepping in the loop on line~\ref{alg:ABAS:ROMLoop}. 
In each iteration of the loop, a time step with the reduced model \eqref{eq:Prelim:ADEIMModel} is taken and the lifted reduced state $\bfV_k\tbfq_k$ is stored.
If the sampling points are to be updated in the current time step $k$, then the predictor model \eqref{eq:Ahead:PredictorModel} is integrated for $C_{\tau}$ time steps for constructing the data sample $\hbfq_k$, the residual matrix $\bfR_k$ is computed, and the sampling points are updated as described in Section~\ref{sec:Prelim:UpdateSamplingPoints}. 
If the sampling points are not updated in the current time step $k$, then the reduced predictor model \eqref{eq:Ahead:ReducedPredictorModel} is integrated for $C_{\tau}$ time steps to compute the data sample $\hbfq_k$ for the data window. 
In line~\ref{alg:ABAS:ApplyBasisUpdate}, the basis update $\bfalpha_k\bfbeta_k^T$ is computed based on the data window $\bfF_k$ and then the adapted basis matrix $\bfV_{k + 1}$ is orthogonalized. 
The empirical-interpolation points $\bfp_k$ to $\bfp_{k + 1}$ are adapted by applying QDEIM\cite{QDEIM} to the adapted basis matrix $\bfV_{k + 1}$. 
Each step in the algorithm scales at most linearly in the dimension $\nh$ of the full-model states.

\section{Numerical experiments}\label{sec:NumExp}
We demonstrate the lookahead strategy on two numerical experiments. 

\begin{figure*}
\begin{tabular}{ccc}
\resizebox{0.33\textwidth}{!}{\huge\input{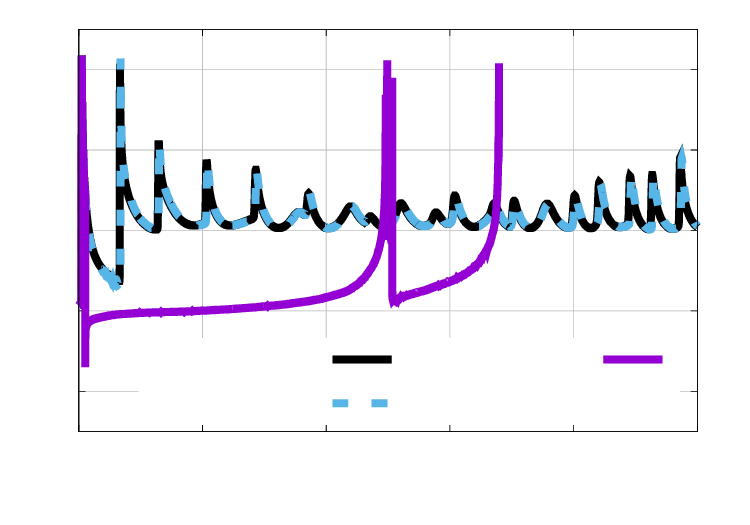}} &
\resizebox{0.33\textwidth}{!}{\huge\input{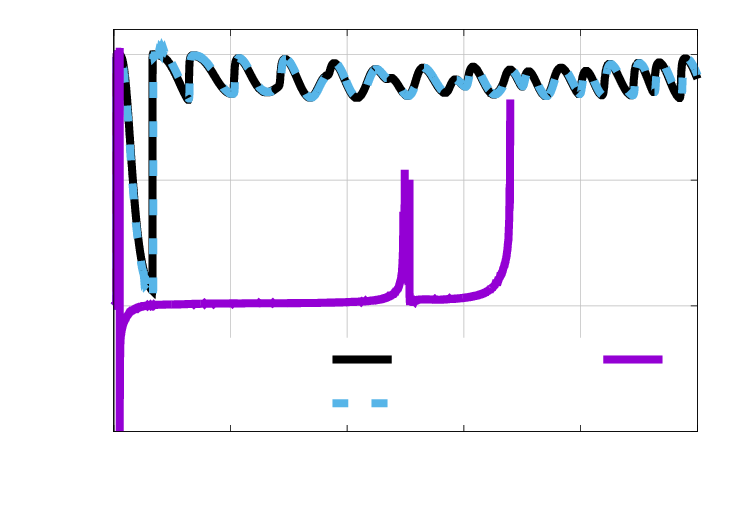}} &
\resizebox{0.33\textwidth}{!}{\huge\input{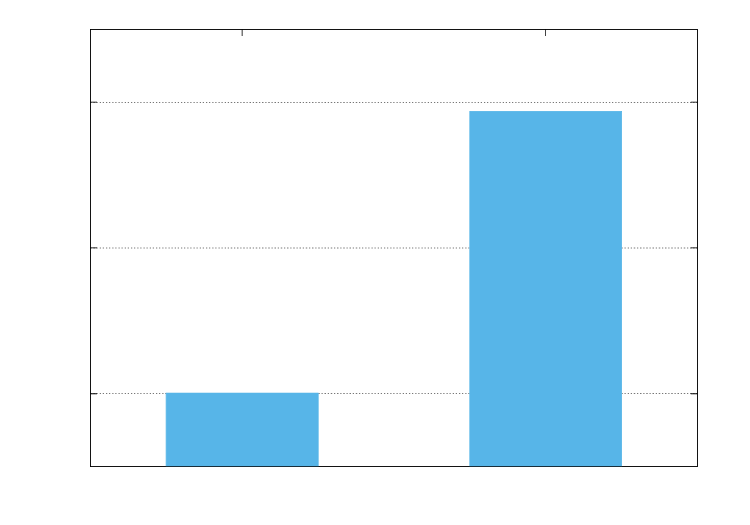}}
\\
(a) intensive property of fluid $\eta$ & (b) combustion progress $\lambda$ &  (c) speedup
\end{tabular}
\caption{Rotating detonation waves: The ADEIM reduced model with lookahead strategy accurately predicts the intensity $\eta$ and progress $\lambda$ at a probe $x = 0.5\pi$ in the spatial domain $\Omega$ and achieves a speedup of two to six in this example. In contrast, ADEIM with the lookback strategy provides a poor approximation.} 
\label{fig:RDE_probe1}
\end{figure*}

\subsection{Rotating detonation waves} \label{sec:NumExp:RDE}
We consider a model for rotating detonation waves \cite{PhysRevE.101.013106} that is motivated by rotating detonation engines. \cite{ANAND2019182,doi:10.1146/annurev-fluid-120720-032612}
In our setting, a single pulse initial condition leads to the formation of a detonation---solution field with sharp gradient---which then spawns a wave. The wave travels over time over a circular domain which then eventually spawns a second wave. The number of waves that are spawned are controlled by the injection parameter $\mu$.

\subsubsection{Governing equations}
The governing equations over the spatial domain $\Omega = [0, 2\pi)$ and time domain $\Tcal = [0, 50]$ are
    \begin{equation} \label{eq:RDEequations}
    \begin{aligned}
\frac{\partial }{\partial t}\eta(x, t)  = & - \eta(x, t) \frac{\partial}{\partial x}\eta(x, t) + \nu \frac{\partial^2 }{\partial x^2}\eta(x, t) \\
 &+ (1- \lambda(x, t)) \omega(\eta(x, t)) + \xi(\eta(x, t))\,,   \\
\frac{\partial }{\partial t}\lambda(x, t)  = &\nu \frac{\partial^2}{\partial x^2}\lambda(x, t) + (1-\lambda(x, t)) \omega(\eta(x, t))\\
& - \beta(\eta(x, t); \mu) \lambda(x, t)\,, 
\end{aligned}
    \end{equation}
    with periodic boundary conditions, spatial coordinate $x \in \Omega$, and time $t \in \Tcal$. The function $\eta: \Omega \times \Tcal \to \mathbb{R}$ represents the intensive property of the working fluid and $\lambda: \Omega \times \Tcal \to [0,1]$  the combustion progress with a value of 1 corresponding to complete combustion. The viscosity parameter is $\nu=0.01$. 
The term $$\omega(\eta(x, t)) = k_{\text{pre}}e^{\frac{\eta(x, t) - \eta_c}{\alpha}}$$ in the model is the heat release function that represents gain depletion. The parameters are set as as follows: the pre-exponential factor is $k_{\text{pre}}=1$, the activation energy is $\alpha=0.3$, and the ignition energy is $\eta_c=1.1$. The injection term is $$\beta(\eta(x, t); \mu) = \frac{\mu}{ 1 + e^{r (\eta(x, t) - \eta_p)}}\,,$$ with $\eta_p=0.5$ and $r=1$. The injection parameter $\mu$ is set to $\mu = 3.5$, if not otherwise noted. The energy loss function is  $\xi(\eta(x, t)) = - \epsilon \eta(x, t)$, with $\epsilon = 0.11$.

\subsubsection{Discretization}
The governing equations \eqref{eq:RDEequations} are discretized in the spatial domain using a first-order upwind scheme at $1024$ equidistant points in the spatial domain $\Omega$, which leads to a model of the form \eqref{eq:contFOM} with states $\bfq(t; \mu)$ of dimension $\nh = 2048$; notice that $\mu$ is the injection parameter that enters the model via the injection term. 
A time-discrete model of the form \eqref{eq:discFOM} is obtained with an implicit Euler discretization with time-step size $\delta t = 10^{-3}$.
The initial condition is $$\eta(x,0) = \frac{3}{2 (1-\cosh (x-1))^{20}}$$ and $\lambda(x,0) = 0 $ for $x \in \Omega$.
Figure~\ref{fig:RDE_FOM2DStills} shows the rotating detonation waves at times $t = 0.05, 0.25, 1, 10, 45$. The dynamics are challenging for model reduction due to the sharp gradients in the solution fields that travel over time as well as the presence of traveling waves that spawn further waves.\cite{Notices}

\subsubsection{Reduced models}
We construct an ADEIM reduced model with state dimension $\nr=9$ and initial window size $\winit = 500$. 
The sampling points are updated every $\nz = 3$ time steps and the number of sampling points are set to $\nrs = \lceil 0.5\nh\rceil$, if not otherwise noted. 
To fill the data window, we use the lookahead strategy described in Section~\ref{sec:Ahead}. For the reduced predictor model \eqref{eq:Ahead:ReducedPredictorModel}, we set $C_{\tau} = 5$, which means that the time-step size $\delta \tau$ of the predictor model is five times smaller than the time-step size $\delta t$ of the full model. The predictor model uses forward Euler as time-integration scheme.  
For comparison, we also show results for an ADEIM model with the lookback strategy following Ref.~\onlinecite{P18AADEIM} and as described in Section~\ref{sec:Prelim:Lookback}; all other parameters of the ADEIM model are the same.

\subsubsection{Results}
Figure~\ref{fig:RDE_probe1} shows a probe of the intensive property of the working fluid $\eta$ and the combustion progress $\lambda$ at location $x = \pi/2$ in the spatial domain $\Omega$. The ADEIM model with the proposed lookahead strategy provides an accurate prediction of the full-model probe whereas the same ADEIM model with the lookback strategy\cite{P18AADEIM} leads to high errors. Plot (c) in Figure~\ref{fig:RDE_probe1} shows the runtime of numerically solving the full model and the ADEIM model with the lookahead strategy. If the full model is discretized on 1024 grid points so that $\nh = 2048$ (two degrees of freedom) then a speedup of about a factor two is achieved with the ADEIM reduced model with lookahead strategy. Discretizing the full model on a fine grid of 2048 grid points so that $\nh = 4096$, a speedup of almost a factor six is obtained. This shows that higher speedups can be expected as the full model becomes more expensive to simulate, which is also in agreement with the results in the following section. 

\begin{figure}[t]
\begin{tabular}{cc}
\resizebox{0.45\columnwidth}{!}{\Huge\input{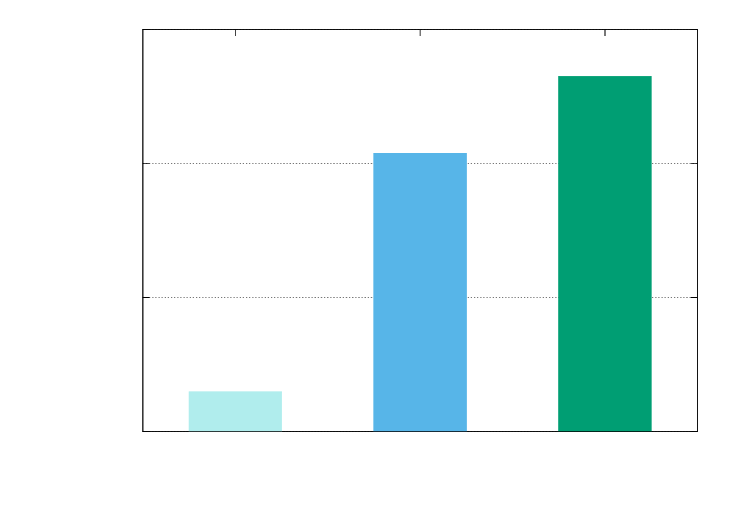}} &
\resizebox{0.45\columnwidth}{!}{\Huge\input{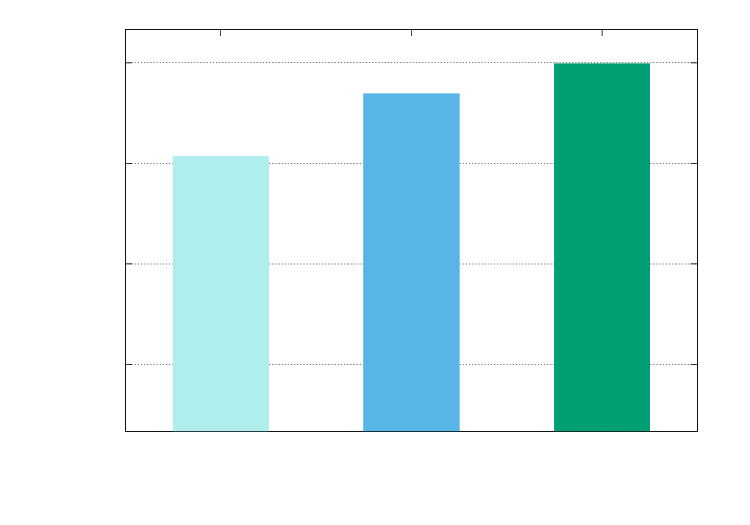}} \\
(a) accuracy & (b) runtime 
\end{tabular}
\caption{Rotating detonation waves: Updating the sampling matrix $\bfS_k$ less frequently (larger $\nz$) leads to lower accuracy in favor of higher speedups.} 
\label{fig:RDE_bar_fixResPct}
\end{figure}

\begin{figure}
\begin{tabular}{cc}
\resizebox{0.45\columnwidth}{!}{\Huge\input{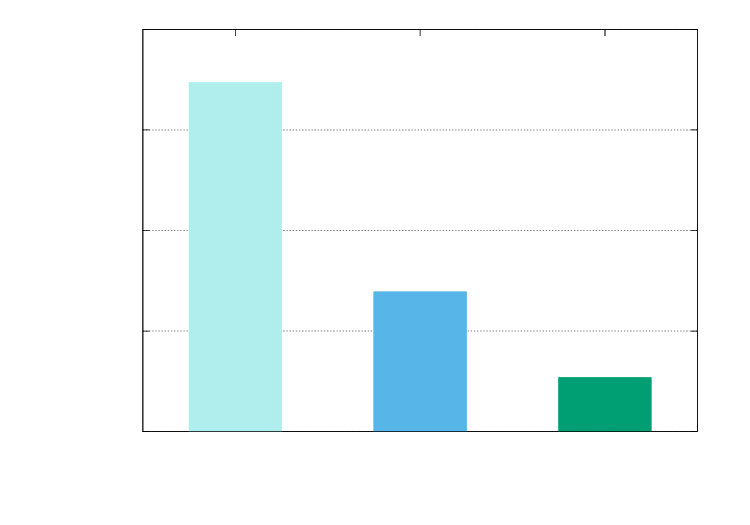}} 
&
\resizebox{0.45\columnwidth}{!}{\Huge\input{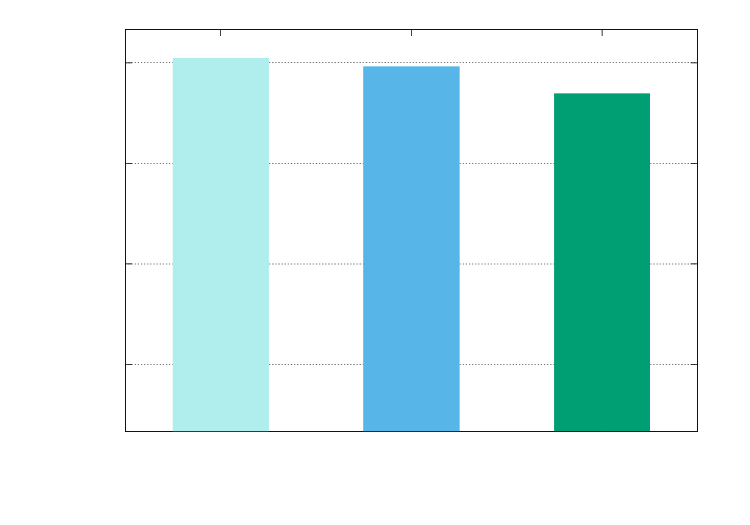}}\\ 
(a) accuracy &  (b) runtime 
\end{tabular}
\caption{Rotating detonation waves: Increasing the number $\nrs$ of sampling points from $\nrs = \lceil 0.3\nh \rceil$ (30\%) to $\nrs = \lceil 0.5\nh \rceil$ (50\%) decreases the error by a factor of about two, which shows that selecting sufficiently many sampling points is key for the accuracy of ADEIM reduced models.} 
\label{fig:RDE_bar_fixUpFreqSamp}
\end{figure}

We now vary the frequency of updating the sampling points $\nz$ and the number of sampling points $\nrs$ and assess the accuracy of the ADEIM reduced model with the average relative error with respect to the full-model states,
    \begin{align}\label{eq:RelErr}
        e(\bfmu) = {\|\tbfQ(\bfmu) - \bfQ(\bfmu)|\|_F^2}/{\|\bfQ(\bfmu)\|_F^2}\,,
    \end{align}
    where $\bfQ(\bfmu)$ is the full model trajectory while $\tbfQ(\bfmu)$ is the ADEIM trajectory.
First, we fix the number of sampling points to $\nrs = \lceil 0.5N\rceil$ and plot the error and the runtime as a function of the frequency of updating the sampling points; see Figure~\ref{fig:RDE_bar_fixResPct}. The plots demonstrate that frequently updating the sampling points results in the most accurate reduced model. However, this comes at the cost of increased simulation time since this implies an increased number of right-hand side evaluations of the full model. 
Second, we fix the frequency of updating the sampling points and vary the proportion of sampling points at which the full model is evaluated. The reduced model errors as well as the speedup with respect to the full model are shown in Figure~\ref{fig:RDE_bar_fixUpFreqSamp}. The plots suggest that for this example, increasing the proportion of sampling points, which reduces the number of solution components that need to be interpolated, increases the accuracy of the reduced model. As $\nrs$ is increased, there is only a slight reduction in speedup.

\subsubsection{Predicting bifurcation diagram}
\begin{figure}
    \resizebox{1.0\columnwidth}{!}{\large\input{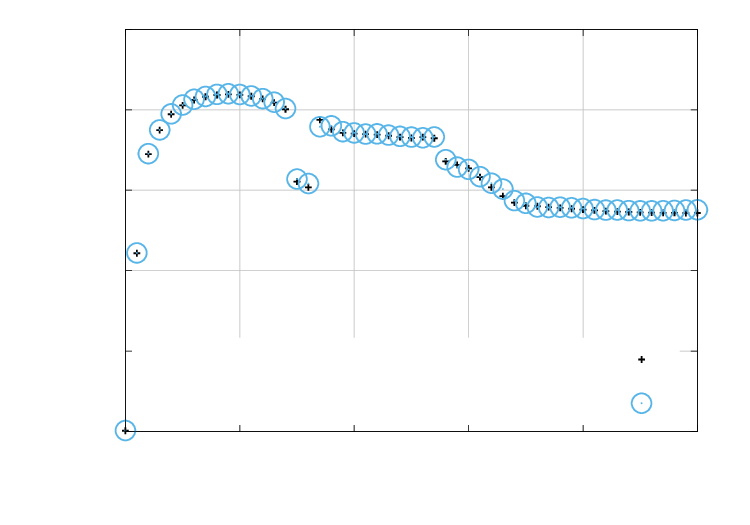}}
\caption{Rotating detonation waves: With the ADEIM reduced model, one can rapidly sweep over the injection parameter $\mu$ to compute the maximum intensity, which indicates how many waves are spawned during the combustion process. The ADEIM model accurately predicts the maximum intensity over a wide range of injection parameters and achieves speedups of about a factor two to six in this example, see Figure~\ref{fig:RDE_probe1}.}
\label{fig:RDE_Bifurcation}
\end{figure}

We now use the ADEIM reduced model with lookahead and the configuration used for generating the probe plots shown in Figure~\ref{fig:RDE_probe1} to predict a bifurction diagram by varying the injection parameter $\mu$. As $\mu$ is varied, the maximum of the intensity $\eta$ over the spatial domain changes, which indicates how many waves are spawned.\cite{ANAND2019182,doi:10.1146/annurev-fluid-120720-032612} Figure~\ref{fig:RDE_Bifurcation} shows that the ADEIM model with lookahead predicts the maximum intensity in close agreement with the full model, while achieving speedups of about a factor two to six in this example.

\begin{figure*}[t]
\begin{tabular}{ccc}
\resizebox{0.33\textwidth}{!}{\huge\includegraphics{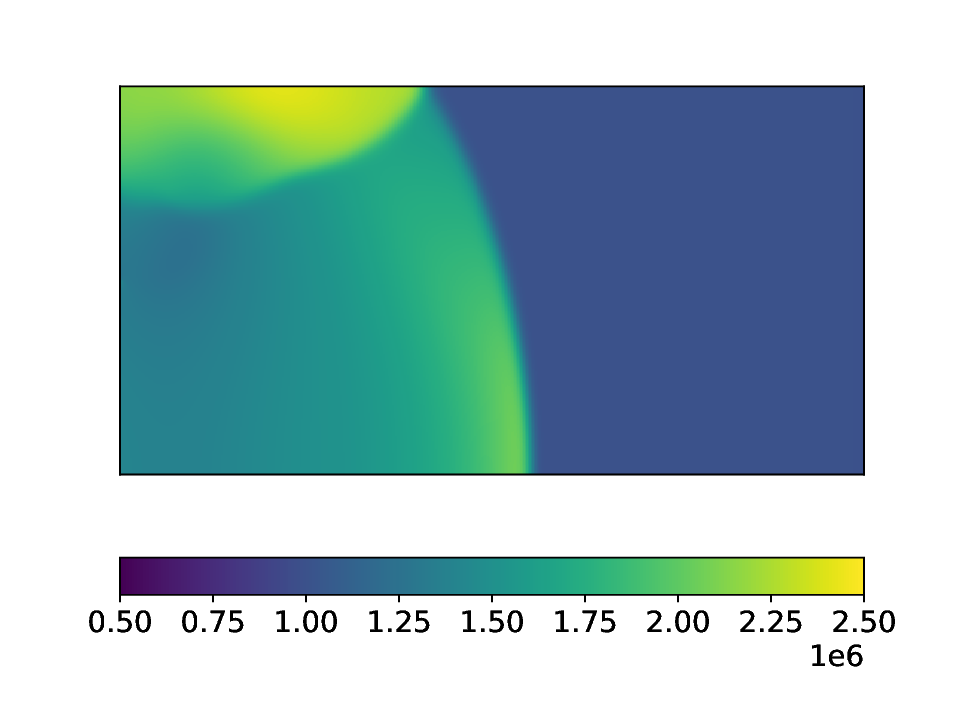}} & \resizebox{0.33\textwidth}{!}{\huge\includegraphics{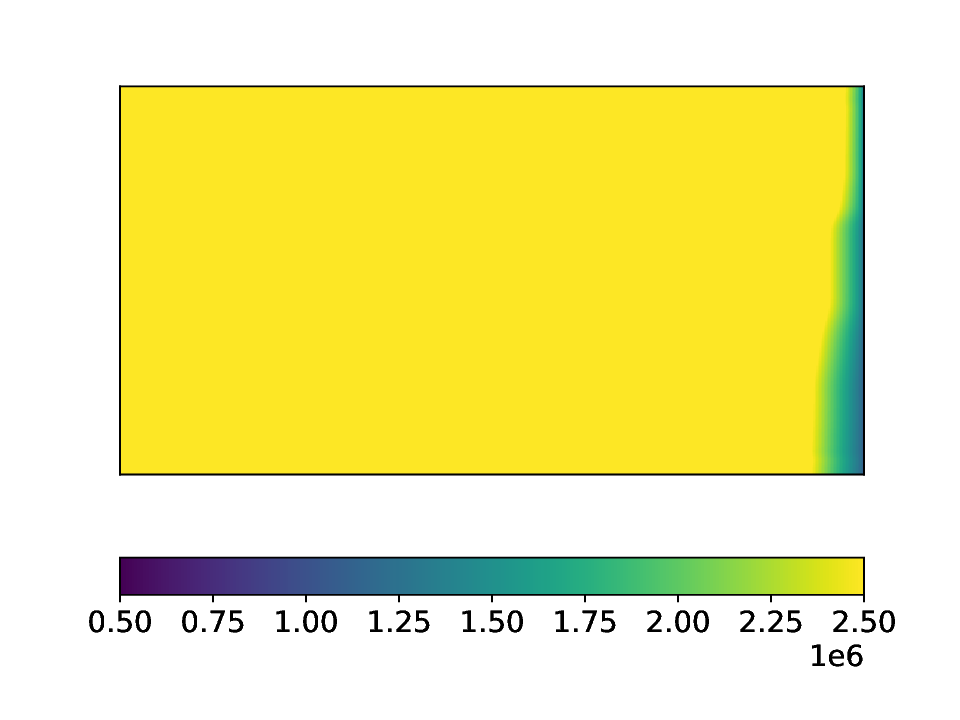}} & \resizebox{0.33\textwidth}{!}{\includegraphics[]{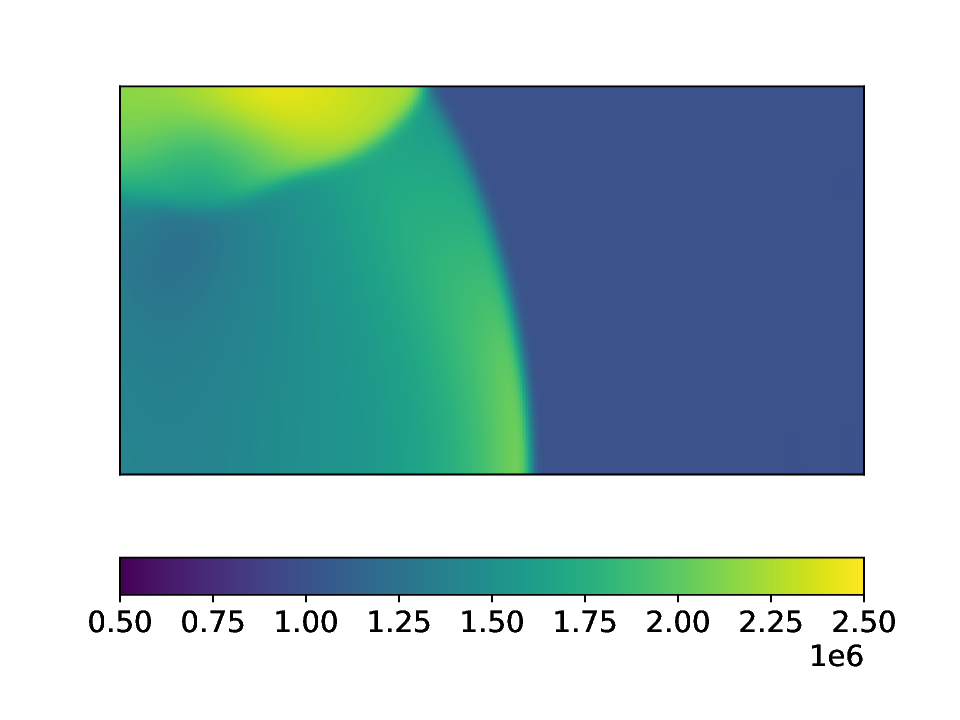}}\\
(a) full model, pressure & (b) ADEIM with lookback, pressure & (c) ADEIM with lookahead, pressure\\
\resizebox{0.33\textwidth}{!}{\includegraphics{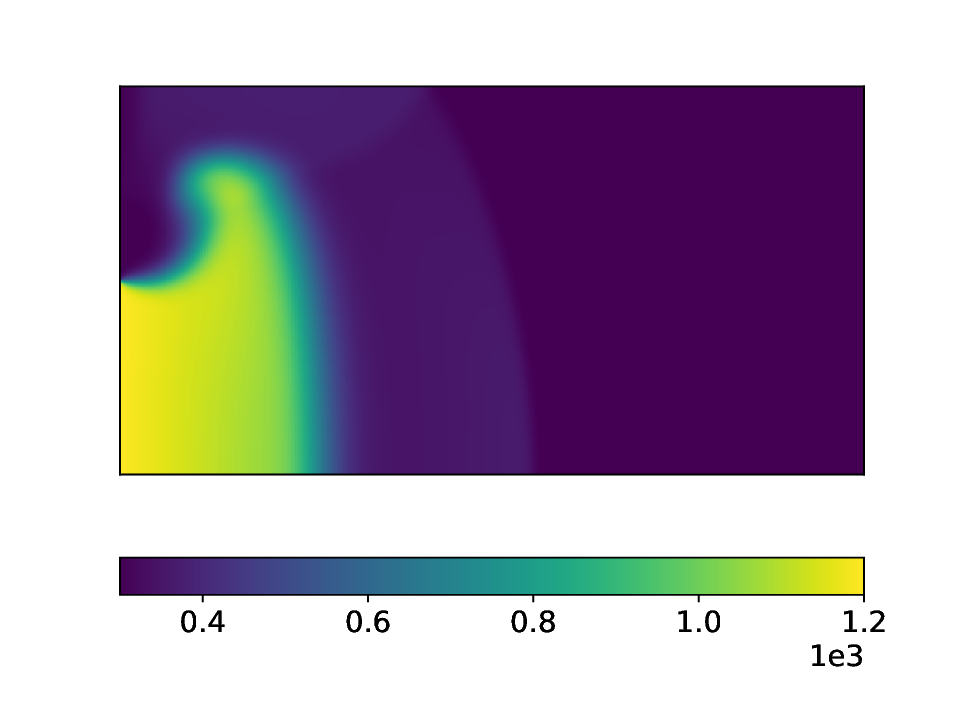}} & \resizebox{0.33\textwidth}{!}{\includegraphics[]{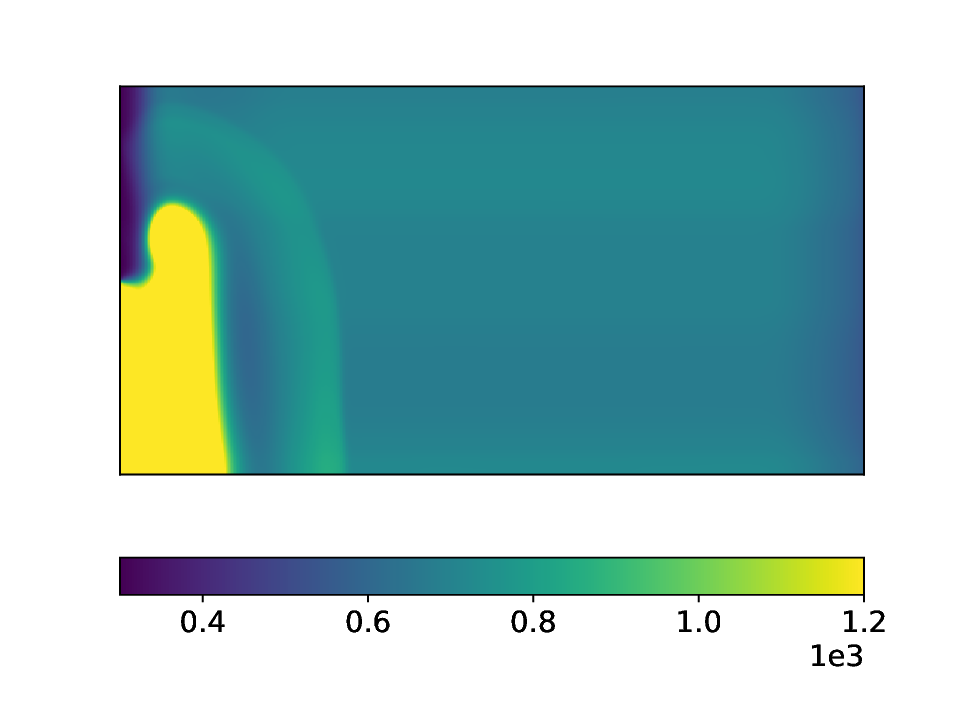}} & \resizebox{0.33\textwidth}{!}{\includegraphics[]{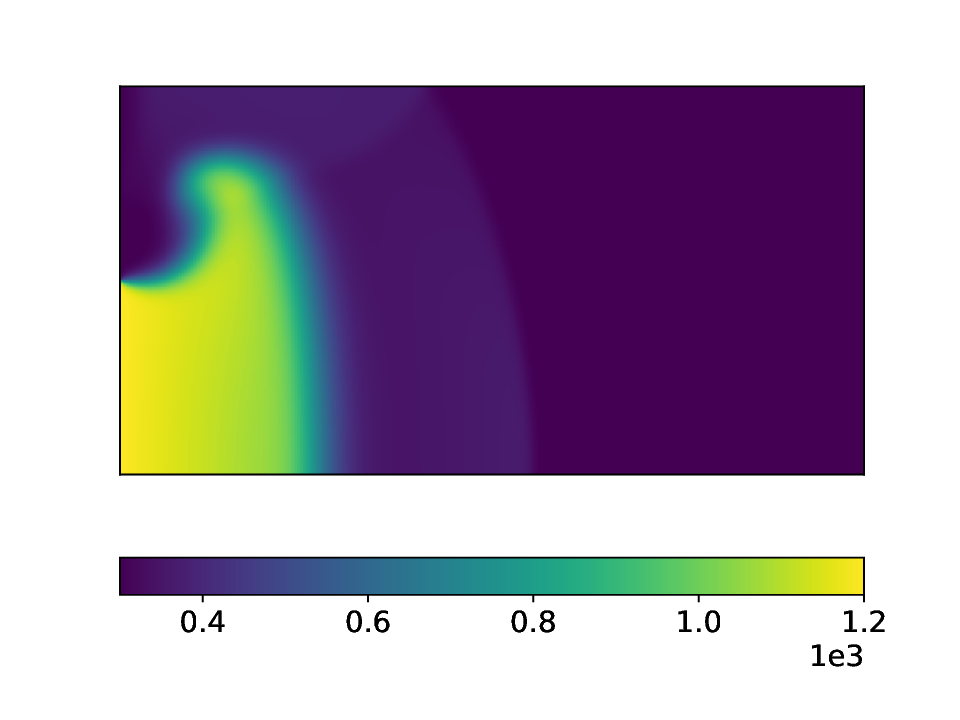}}\\
(d) full model, temperature & (e) ADEIM with lookback, temperature & (f) ADEIM with lookahead, temperature
\end{tabular}
\caption{Mixing layer: The ADEIM reduced model with the proposed lookahead strategy accurately approximates the pressure and temperature fields of the full model, whereas ADEIM with lookback leads to poor approximations and eventually becomes unstable.} 
\label{fig:MixingLayer_rhseval_instability}
\end{figure*}

\begin{figure*}[t]
\begin{tabular}{c}
\resizebox{1.0\textwidth}{!}{\large\input{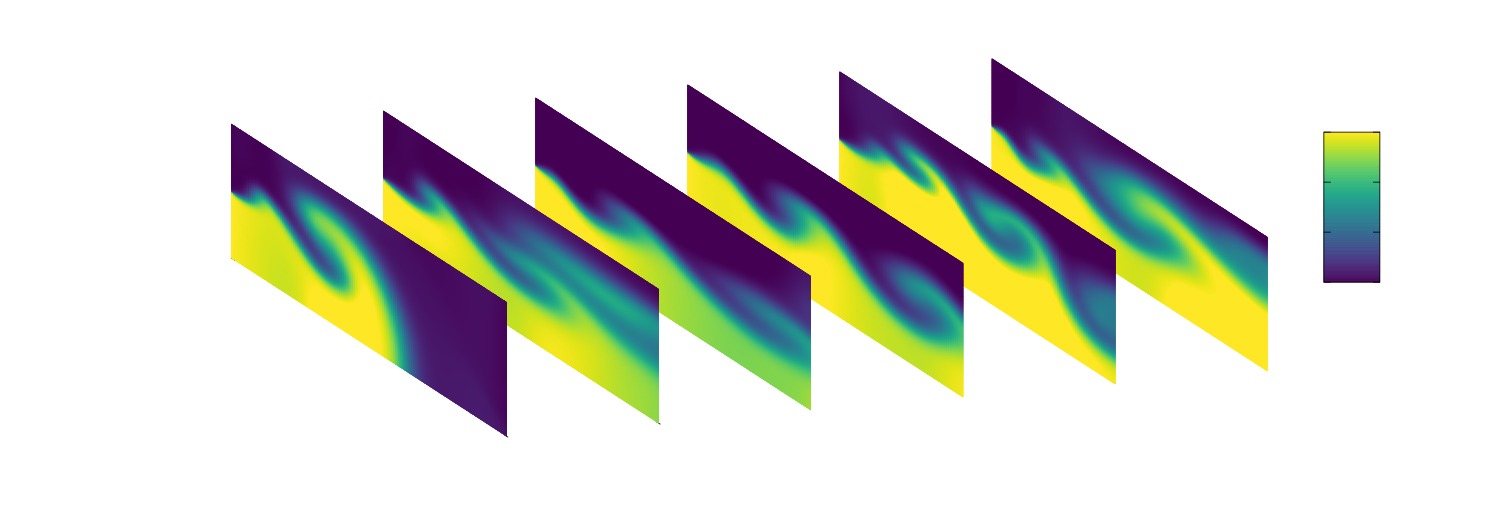}}\\
(a) full model\\
\resizebox{1.0\textwidth}{!}{\large\input{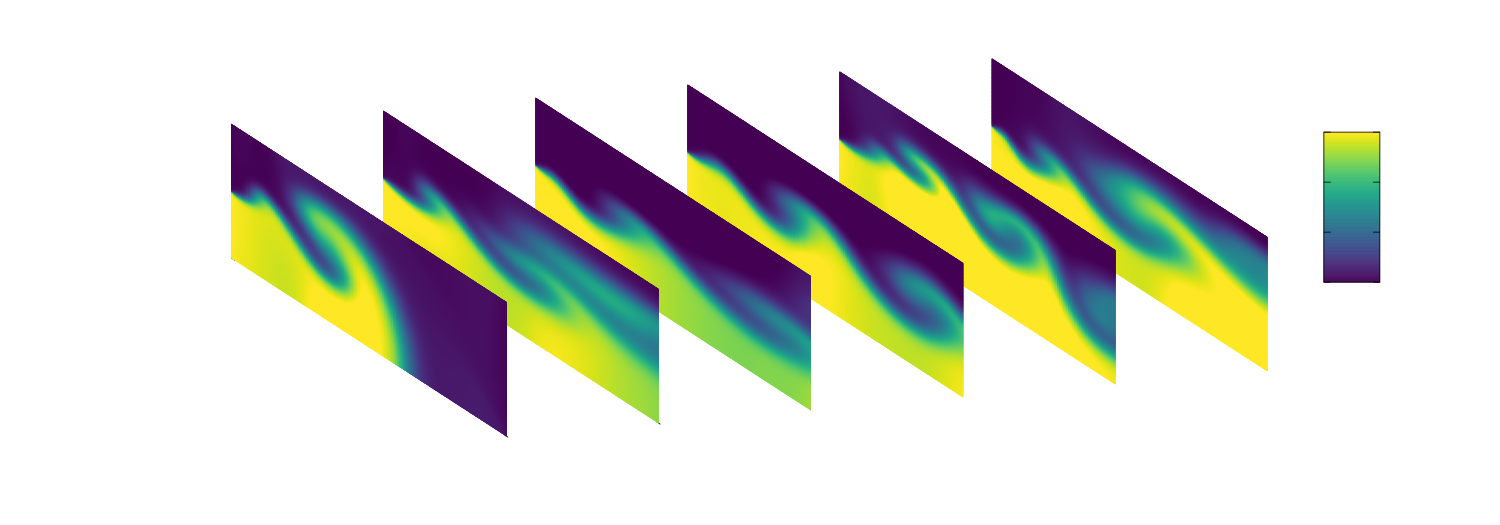}}\\
(b) ADEIM with lookahead strategy predictions
\end{tabular}
\caption{Mixing layer: The ADEIM reduced model accurately predicts the temperature field of the full model.}
\label{fig:Mixing:3DTemperature}
\end{figure*}

\begin{figure*}[t]
\begin{tabular}{c}
\resizebox{1.0\textwidth}{!}{\Huge\input{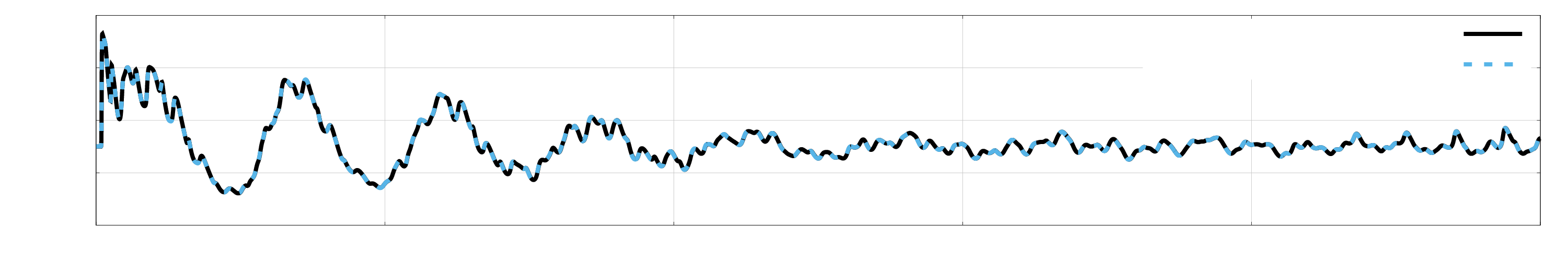}}\\
(a) pressure probe, mixing layer with low temperature (case 1)\\
\resizebox{1.0\textwidth}{!}{\Huge\input{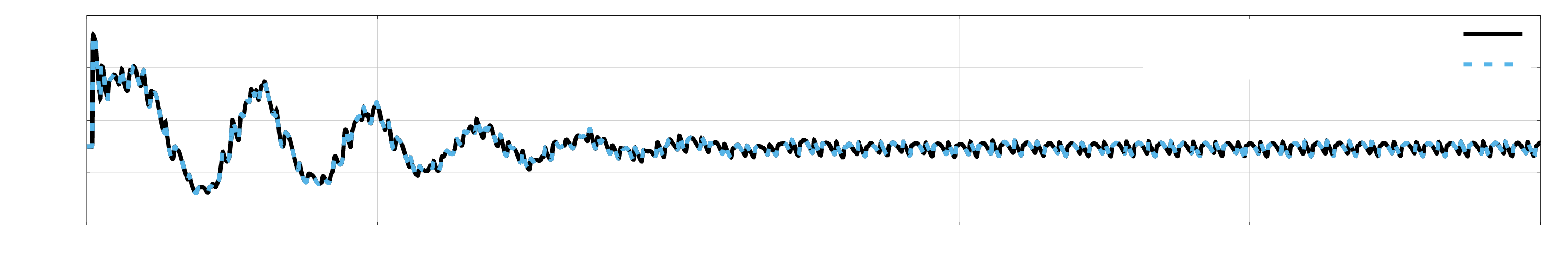}}\\
(b) pressure probe, mixing layer with high temperature (case 2)
\end{tabular}
\caption{Mixing layer: The ADEIM reduced model with the proposed lookahead strategy accurately approximates the pressure at the probe location over the time window $[0, 0.01]$.}
\label{fig:NumExp:Mixing:Probes}
\end{figure*}

\subsection{Mixing layer flow}
We consider a fluid-flow model that describes a mixing layer, which is of interest in understanding turbulence and combustion in aerospace engineering and beyond.\cite{Lesieur1995}
Two spatially separated inlets inject fluid at different temperatures and velocities, which causes the pressure to oscillate. The oscillations are controlled by the ratio of temperature and velocity between the two inlet jets.

\subsubsection{Setup}
In our setup, the spatial domain is of 0.1m length and 0.04m height so that $\Omega = [0.1] \times [0.04] \subset \mathbb{R}^2$. The time domain is $[0, 0.01]$. 
The flow is governed by an equation that we write in conservative form as 
\begin{equation}\label{eq:NumExp:NSPDE}
\partial_t \vec{q}(t, x, y) + \nabla \cdot \left(F(t, x, y, \vec{q}) - F_v(t, x, y, \vec{q})\right) = 0 
\end{equation}
with $\vec{q} = [\rho;  \rho v_x; \rho v_y; \rho e]$, where $\rho$ is the density, $v_x$ and $v_y$ are the velocity in $x$ and $y$ directions, respectively, and $e$ is the total energy. The inviscid flux $F$ is
\[
F = \begin{bmatrix}
\rho v_x\\
\rho v_x^2 + p\\
\rho v_x v_y\\
\rho v_xe + pv_x\\
\end{bmatrix}\vec{i} + \begin{bmatrix}
\rho v_y\\
\rho v_x v_y\\
\rho v_y^2 + p\\
\rho v_ye + pv_y\\
\end{bmatrix}\vec{j}\,,
\]
where $p$ is the pressure. The total energy and the pressure are related as \begin{equation}\label{eq:NumExp:EnergyEOS}
p = (\gamma - 1) (\rho e - 0.5\rho(v_x^2 + v_y^2)),
\end{equation} where $\gamma = c_p/(c_p - R)$ is the specific heat ratio, $c_p = 1.538$ [kJ/kg/K] is the specific heat capacity of the fluid, and $R = R_u/M_w$. The quantity $R_u = 8.314$ [J/(mol $\cdot$ K)] is the universal gas constant and the molecular weight of the fluid is set to $M_w = 21.32$ [g/mol]. 
The viscous flux $F_v$ is
\[
F_v = \begin{bmatrix}0\\
\tau_{xx}\\
\tau_{xy}\\
\tau_{xx}v_x + \tau_{xy}v_y - j_x^q\end{bmatrix}\vec{i} + \begin{bmatrix}0\\
\tau_{xy}\\
\tau_{yy}\\
\tau_{xy}v_x + \tau_{yy}v_y - j_y^q\end{bmatrix}\vec{j}\,,
\]
with the viscous shear tensor
\[ 
\tau = \begin{bmatrix}\tau_{xx} & \tau_{xy}\\ \tau_{xy} & \tau_{yy}\\\end{bmatrix} =
\eta \begin{bmatrix} \frac{1}{3} \frac{\partial v_x}{\partial x} & \frac{\partial v_x}{\partial y} + \frac{\partial v_y}{\partial x}\\\frac{\partial v_x}{\partial y} + \frac{\partial v_y}{\partial x} & \frac{1}{3} \frac{\partial v_y}{\partial y}\end{bmatrix}
\]
and the mixture viscosity coefficient $\eta$, which we set to $\eta = 7.35 \times 10^{-4}$. The vector $j^q = [j_x^q; j_y^q]$ is the diffusive heat flux vector defined as
\[
j^q = - \kp \nabla T = -\kp \begin{bmatrix}\frac{\partial T}{\partial x} \\ \frac{\partial T}{\partial y}\end{bmatrix}\,,
\]
which depends on the thermal conductivity $\kp$ and the temperature $T$. The temperature $T$ satisfies the relationship 
\begin{equation}\label{eq:NumExp:PressTempRho}
\rho = \frac{p M_w}{R_u T}\,.
\end{equation}
The thermal conductivity is $\kp = \mu_p c_p/P_r$, where $\mu_p = 7.35 \times 10^{-4}$ and $P_r = 0.713$ is a dimensionless quantity. 

On the left boundary, we impose inlet boundary conditions. The y-component of the velocity is set to $v_y = 0$ on the left boundary. In the top half of the domain from $y = 0.02$ to $y = 0.04$, we set the x-component of the velocity to $v_x = 100$ [m/s] and the temperature to $T = 300$ [K]. In the bottom half of the domain from $y = 0$ to $y = 0.02$, the x-component of the velocity is set to $v_x = 400$ [m/s]. For the temperature in the bottom half of the domain, we distinguish between two cases. In case one, the low-temperature case, we set the temperature to $T = 700$ [K]. In case two, the high-temperature case, we set the temperature to $T = 1200$ [K]. 
In both cases, the pressure $p$ at the inlet is zeroth order extrapolated from the interior by a homogeneous Neumann boundary condition with $\partial p/\partial \vec{n} = 0$, where $\vec{n}$ is the normal vector to the left boundary. 
Density and energy at the inlet are calculated based on the temperature, pressure, energy relationship given in \eqref{eq:NumExp:EnergyEOS} and \eqref{eq:NumExp:PressTempRho}.

On the right outlet boundary, we impose homogeneous Neumann boundary condition $\partial v_x/\partial \vec{n} = 0$, $\partial v_y/\partial \vec{n} = 0$, $\partial T/\partial \vec{n} = 0$, where $\vec{n}$ is the normal to the right boundary. We set 
the pressure to constant one at the outlet. Density and energy are calculated from \eqref{eq:NumExp:EnergyEOS} and \eqref{eq:NumExp:PressTempRho}. On the top and bottom boundaries, we impose slip wall conditions. 

The initial condition throughout the domain sets $p = 1$ [MPa], $T = 300$ [K], $v_x = 100$ [m/s], and $v_y = 0$ [m/s]. Density and  energy are calculated from \eqref{eq:NumExp:EnergyEOS} and \eqref{eq:NumExp:PressTempRho}.

We discretize the equations \eqref{eq:NumExp:NSPDE}
with a Rusanov semi-discrete scheme, which is of first order. The spatial domain is discretized with a $500 \times 500$ grid, with $\Delta x = 0.0002$ [m] and $\Delta y = 0.00008$ [m]. Thus, the state dimension is $\nh = 1000000$ over all four conserved variables combined. The time integrator is TR-BDF2, which is a second-order implicit method, with $\Delta t = 10^{-8}$ [s].

\subsubsection{Reduced models}
We construct an ADEIM reduced model with state dimension $\nr = 6$ and initial window size $\winit = 100$. The sampling points are updated every $\nz = 4$ time steps and the number of sampling points are set to $\nrs = \lceil 0.001 \nh \rceil$, which means that $0.1\%$ of the component functions of the full-model right-hand side function are sampled for adapting the basis. 
We consider an ADEIM model with the lookback strategy as described in Section~\ref{sec:Prelim:ADEIM} and an ADEIM model with the lookahead strategy of Section~\ref{sec:Ahead}. For the lookahead strategy, we set $C_{\tau} = 5$, if not stated otherwise. The predictor model uses forward Euler time integration.

\begin{figure}[t]
\resizebox{0.8\columnwidth}{!}{\LARGE\input{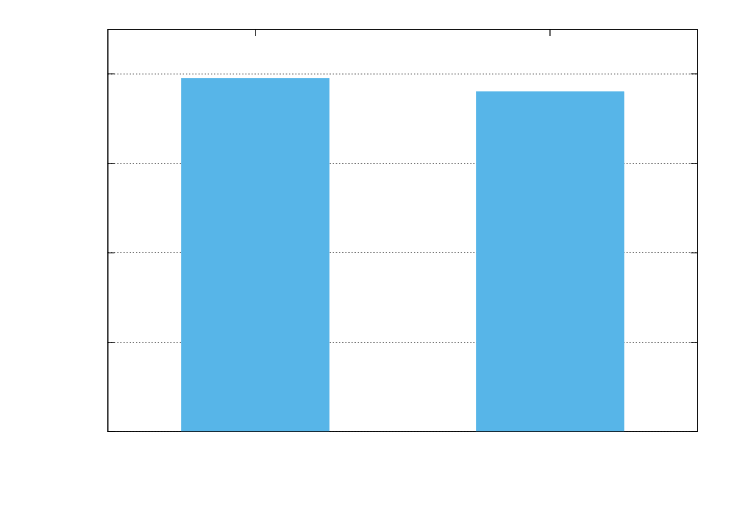}}
\caption{Mixing layer: The ADEIM reduced model with the proposed lookahead strategy achieves a speedup of almost 14 for the high- and low-temperature case of the mixing layer setup.}
\label{fig:MixingLayer:Speedup}
\end{figure}

\subsubsection{Results}
\begin{figure}
\begin{tabular}{cc}
\resizebox{0.48\columnwidth}{!}{\Huge\input{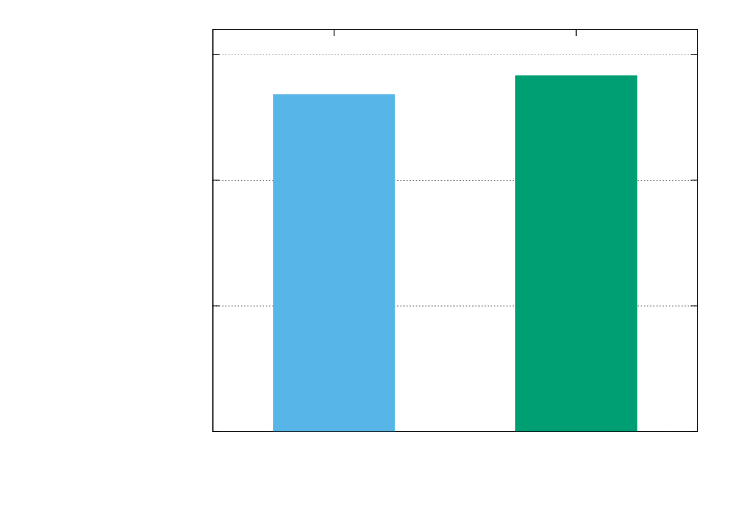}} & 
\resizebox{0.48\columnwidth}{!}{\Huge\input{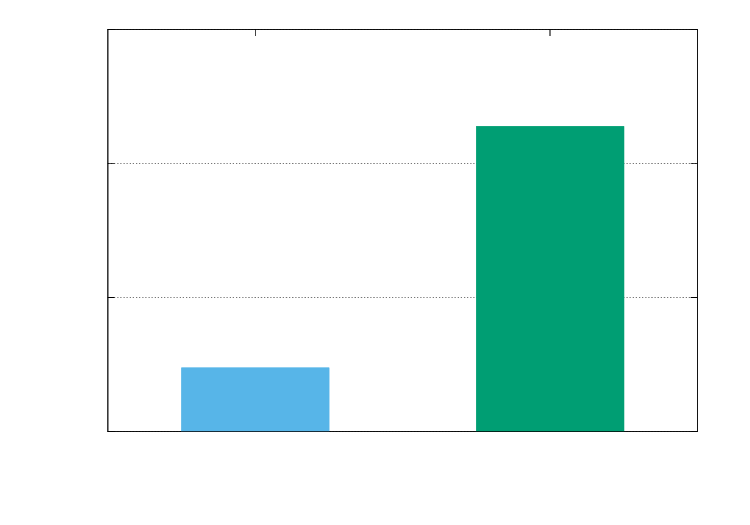}}\\
(a) error & (b) speedup
\end{tabular}
\caption{Mixing layer: Increasing the time-step size $\delta \tau$ of the predictor model has little effect on the error of the ADEIM reduced model (see plot a) but increases the speedup to more than 20. The results show that a crude predictor model is sufficient to guide the adaptation of the ADEIM reduced space in this example.}
\label{fig:MixingLayer:VaryingDt}
\end{figure}

Figure~\ref{fig:MixingLayer_rhseval_instability} shows the pressure and temperature fields at time $t = 0.00008$ [s] for the high-temperature case computed with the full model, the ADEIM reduced model with lookback, and the ADEIM reduced model with lookahead. The results show that the lookahead strategy leads to an accurate approximation of the full-model pressure and temperature fields, whereas the reduced models with lookback lead to poor approximations and eventually become unstable. Figure~\ref{fig:Mixing:3DTemperature} shows the temperature field at several other points in time. The ADEIM reduced model accurately predicts the full-model fields. 
Figure~\ref{fig:NumExp:Mixing:Probes} shows the pressure at probe location $x = 0.025$ and $y = 0.02$ in the spatial domain $\Omega$ for the low- and high-temperature case. The ADEIM model with lookahead accurately predicts the pressure at the probe location for both configurations, even though the two configurations lead to distinctly different pressure oscillations.
The speedup of the ADEIM reduced model with lookahead strategy compared to the full model is shown in Figure~\ref{fig:MixingLayer:Speedup}. A speedup of almost 14 is achieved with the ADEIM reduced model.
In Figure~\ref{fig:MixingLayer:VaryingDt}, we show that decreasing $C_{\tau}$ from $C_{\tau} = 5$ to $C_{\tau} = 2$, so that the time-step size of the predictor model $\delta \tau$ increases, barely changes the relative average error \eqref{eq:RelErr} while increasing the speedup to more than 20. The results show that a crude predictor model is sufficient to guide the adaptation of the ADEIM reduced space in this example.

\section{Conclusions}\label{sec:Conc}
The data that are gathered from the full model have a major impact on the quality of online adaptive reduced models because the data determine towards which dynamics the reduced spaces are adapted. The lookahead strategies proposed in this work aim to predict the dynamics that are likely to be seen in the immediate future so that the full model can be queried for informative data. The numerical experiments show that the online adaptive reduced models with lookahead strategies are more accurate and stable, even when online adaptive reduced models with previously introduced data-gathering strategies fail to be predictive.

\section*{Acknowledgements} The authors would like to thank Cheng Huang and Christopher Wentland for guidance on setting up the mixing layer numerical experiment.  This work was supported in part by the Air Force Center of Excellence on Multi-Fidelity Modeling of Rocket Combustor Dynamics, award FA9550-17-1-0195. B.P.~was additionally partially supported by the Air Force Office of Scientific Research (AFOSR) award FA9550-21-1-0222 (Dr.~Fariba Fahroo). This work was also supported in part through the NYU IT High Performance Computing resources, services, and staff expertise.

\clearpage
\bibliographystyle{numeric}
\bibliography{main.bib}

\end{document}